\documentclass[12pt]{article}
\usepackage[dvipdfmx]{graphicx} % Required for inserting images
\usepackage[utf8]{inputenc}
\usepackage{amsthm}
\usepackage{amsmath}
\usepackage{natbib}
\usepackage{color}
\usepackage{ulem}
\usepackage{url}
\usepackage{bm}
\usepackage{here}
\usepackage{float}

\begin{document}
\title{Elongation of material lines and vortices by Euler flows on two-dimensional Riemannian manifolds}
\author{Koki Ryono and Keiichi Ishioka
\\Graduate School of Science, Kyoto University, \\
Kitashirakawa-Oiwake-cho, Sakyo-ku, Kyoto 606-8502, Japan}
\date{March 2025}

\begin{abstract}
We study the influence of the differential geometry of the flow domain on the motion of fluids on two-dimensional Riemannian manifolds, particularly on the elongation of material lines and vortices. We derive a formula for the second order time derivative of the square of the distance between close fluid particles and show that a curvature term appears. The elongation of a material line is accelerated by negative curvature. The use of this expression extends Haller's definition of hyperbolic domains to flows on curved surfaces. The need to consider curvature effects is illustrated by three examples. The example of a curved two-dimensional torus implies that the filamentation of vortices can be triggered by negative curvature.
\end{abstract}

\maketitle

\, \\
\textit{This is the version of the article before peer review or editing, as submitted by authors to Fluid Dynamics Research. IOP Publishing Ltd is not responsible for any errors or omissions in this version of the manuscript or any version derived from it.}

\section{Introduction}
Although often overlooked, fluid motion is strongly influenced by the geometric features of the flow domain. These geometric features can be broadly divided into three main categories: topology, symmetry, and metrics. The topology of a flow domain imposes restrictions at the most fundamental level on what kind of flow vector fields are possible. Also, in compact flow domains, the elliptic operators linking the stream function and vorticity have discrete spectra, which is one reason why the low wavenumber mode structures are important, as seen in the large-scale motions of the Earth's atmosphere. In the Earth's oceans, Kelvin waves traveling along the coast are due in part to the Earth's rotation, but also to the existence of boundaries in the flow domain. Second, the symmetry of the flow domain provides a constraint; Noether's theorem guarantees the existence of conserved quantities. For flows on a sphere, the three components of angular momentum are conserved, corresponding to rotational symmetry in the three directions.  The conservation of angular momentum is an effect that cannot be overlooked in the formation of large jets in the Earth's atmosphere. The third category, metrics or differential geometrical property, is perhaps less often considered than the two types above, but Rossby waves in geophysical fluids are waves generated by the latitudinal variation of the vorticity of the solid rotating flow that represents the Earth's rotation. This is an effect of the curvature of the Earth and can be seen as a manifestation of the metric effect.

In this study, we explore how the differential geometrical features of the flow domain fundamentally influence the motion of the fluid. Specifically, we examine the effect of the curvature of a Riemannian manifold on fluid motion.
As is well known in differential geometry, the behavior of two adjacent geodesics is determined by the curvature of the manifolds \citep[see e.g.,][]{arnold2009topological, milnor1963morse}. That is, in a manifold of positive curvature, the two geodesics will eventually converge, while in a manifold of negative curvature, they will eventually diverge. This is described by the Jacobi equation for the variational vector. Therefore, on a manifold of negative curvature, the motion of a mass point that is not subject to external forces is highly dependent on its initial velocity and highly uncertain. We adopt the Euler equation as the simplest model of a fluid on a Riemannian manifold. The Euler equation can be considered as a modification of the geodesic equation, in which the cluster of mass points is made continuous and the pressure force is introduced to satisfy incompressibility \citep{mitsumatsu2008}. Therefore, the Euler fluid on Riemannian manifolds is expected to be under a similar effect to that described by the Jacobi equation. In other words, the acceleration of the distance between particles in the fluid according to the Euler equation is expected to be affected by the curvature of the manifold.

The question of how the distance between nearby particles changes with time and, or almost equivalently, the rate at which a material line embedded in a fluid stretches and contracts, is particularly important in two-dimensional fluids. In two-dimensional Euler fluids, the vorticity is frozen to the particles. The filamentation process, in which the vortices are deformed and stretched into an elongated shape, is an important process for two-dimensional flows to become turbulent. As small-scale vortices in the vorticity field are elongated, enstrophy is transported in the high wavenumber direction \citep[see e.g.,][]{vallis2017atmospheric}, while energy is transported in the low wavenumber direction. In addition, the more complex mixing of the filamentary vortices leads to the formation of macroscopic structures. Therefore, it is important to study the behavior of distances between fluid particles in close proximity in order to quantify the filamentation and mixing due to vorticity advection. In other words, filamentation and mixing of the vorticity field occurs mainly in regions where the distance between fluid particles increases at an accelerating rate. Determining in which regions of the flow field such filamentation is active, or conversely, in which regions the vortices retain their identity, has been well studied. Regions of the former type are called hyperbolic domains, while those of the latter type are called elliptical or vortex domains, and various definitions have been proposed. The Okubo-Weiss parameter \citep{okubo1970horizontal,weiss1991dynamics}, which compares the magnitude of strain velocity and vorticity, is a typical example, which is calculated diagnostically from Eulerian field quantities such as velocity and vorticity. In this sense, hyperbolic domains determined by the Okubo-Weiss parameter are classified as Eulerian definitions of hyperbolic domains. On the other hand, the determination of hyperbolic domains using the finite-time Lyapunov exponent and Haller's definition of hyperbolic domains, which will be discussed below, are based on the motion of particles and correspond to Lagrangian definitions of hyperbolic domains.

In the present study, we derive an expression for the second order time derivative of the square of the distance between close fluid particles for an Euler fluid on a Riemannian manifold and show that the curvature of the manifold affects the time change of the distance between fluid particles. The derived expression is a quadratic form corresponding to the strain acceleration tensor. The same type of formula has been used to determine the Lagrangian hyperbolic domain of \cite{haller2001lagrangian} and \cite{haller2005objective}, in the case of flat domains. In these studies, Haller proposed a method for determining hyperbolic and elliptical (vortex) domains of a flow field based on the variation of the distance between fluid particles in close proximity. Therefore, the formula we derive below extends the hyperbolic domain proposed by Haller to Euler fluids on curved surfaces.
In particular, the determination of hyperbolic domains on spherical surfaces has important implications for applications such as the advection of chemical tracers on the Earth. Furthermore, the equations we derive will motivate the study of Euler flows on manifolds with regions of negative curvature. There have been several studies of vortices on negatively curved surfaces from the point of view of superfluid physics \citep[][etc.]{turner2010vortices, reuther2015interplay}. The main interest in these studies was to investigate the relationship between the curvature and the potential energy produced by point vortices, and to study the topology of the streamlines of steady-state solutions of the Navier-Stokes equations. However, our research is motivated by an interest in the unsteady motion of continuous vortices of finite size. Namely, in regions of negative curvature, the elongation of the material line advected by the Eulerian flow can be accelerated by the curvature effect, thereby promoting strong mixing.

In two-dimensional fluids on planar and spherical surfaces, it is well known that coherent vortices form after long periods of turbulence. On a flat torus or a sphere, coherent vortices that survive have strong radial gradients of vorticity, called mixing barriers \citep[e.g.,][]{dritschel2015late}. The question of whether vortices that survive longer on manifolds of negative curvature would have a different structure than on flat or spherical surfaces is interesting in the context of pattern formation of two-dimensional flows. Although this question cannot be answered by the present study, the equations we derive below will provide a tool to explore such a question.

This paper is organized as follows. First, in section 2, we write down the Euler equation on Riemannian manifolds and give two special concrete examples, the standard sphere and a curved torus. We also review (section 2.2) the definition of hyperbolic type domains in the plane as defined by Haller. In section 2.3, we derive expressions for the time variation of the elongation of material lines on general Riemannian manifolds by differential geometrical calculations (section 2.3.1) and extend Haller's definition of hyperbolic type domains (section 2.3.2). In section 3, we discuss the necessity and usefulness of the equations derived in section 2.3 using three specific examples of flow fields: The first example (section 3.1) is a simple jet on a sphere to verify that the formula derived in section 2.3 is consistent with the explicit expression for the length of the advected material line. The second example (section 3.2) is a large-scale vortex on a sphere, where we see that Haller's hyperbolic domain changes with and without considering curvature effects. The third example (section 3.3) is the onset of turbulence on a curved torus. In this example, we see the process of filamentation of the vortices initially present in the region of negative curvature. Using the equation derived in section 2.3, we show that the filamentation is influenced by the negative curvature. Finally, in section 4, we give a discussion based on the previous sections and a conclusion to the whole paper.

\section{Theory}
In this section, we begin by introducing the Euler equation on Riemannian manifolds and attempt to extend Haller's definition of hyperbolic type domains on Riemannian manifolds.
Note that the notation used below to describe Riemannian geometry follows textbooks such as
\cite{nishikawa2002geometric}.

\subsection{Euler equation on a Riemannian manifold}
Let \(M\) be a smooth Riemannian manifold, endowed with the Riemannian metric \(g=(g_{ij})\). The Levi-Civita connection \(\nabla\) is uniquely determined from \(g\), and the Euler equation on \(M\) is given by:
\begin{align}
    \frac{\partial \bm{u}}{\partial t} + \nabla_{\bm{u}} \bm{u} = -{\rm grad}\,p,\quad {\rm div} \,\bm{u} = 0\label{eulereq}
\end{align}
where \(t\) is the time, \(\bm{u}\) is the vector field of the flow, and \(p\) is the pressure field
\citep{arnold2009topological}. Recall that the \(\nabla_X Y\) for two vector fields \(X\) and \(Y\) is the covariant derivative of \(Y\) in the direction of \(X\) defined by connection \(\nabla\). Let \(P\in M\) be an arbitrary point, and take an open neighborhood \(U\) of \(P\), and suppose that \(x=(x^1, \cdots, x^n)\) be a coordinate system on \(U\), then we can use a natural frame:
\begin{align*}
    \partial_1 = \left(\frac{\partial}{\partial x^1}\right)_P, \cdots, \partial_n = \left(\frac{\partial}{\partial x^n}\right)_P
\end{align*}
as the basis of the tangent space \(T_P M\). In \eqref{eulereq}, the connection \(\nabla\) is the Levi-Civita connection, which can be explicitly written as \(\nabla_{\partial i}\partial_j = \sum_{k=1}^n \Gamma_{ij}^k \partial_k\), by using the connection coefficients:
\begin{align*}
    \Gamma_{ij}^k = \frac{1}{2} \sum_{l=1}^n g^{kl}\left(\frac{\partial g_{il}}{\partial x^j}+ \frac{\partial g_{jl}}{\partial x^i}- \frac{\partial g_{ij}}{\partial x^l}\right),
\end{align*}
where we denote the inverse of \((g_{ij})\) by \((g^{ij})\).  The Euler equation \eqref{eulereq} expressed in components is written as:
\begin{align}
    \frac{\partial u^k}{\partial t} + \sum_{i=1}^n u^i\frac{\partial u^k}{\partial x^i} + \sum_{i=1}^n\sum_{j=1}^n \Gamma_{ij}^k u^i u^j = - \sum_{i=1}^n g^{ik}\frac{\partial p}{\partial x^i}\quad(k=1,\cdots, n).\label{eulereq_component}
\end{align}
If we consider a two-dimensional manifold, as we will do in most part of the present study, by taking the curl of \eqref{eulereq_component} we obtain the following vorticity equation:
\begin{align}
    \frac{Dq}{Dt} := \frac{\partial q}{\partial t} + \langle {\rm grad}\,q, u\rangle = 0.\label{vorticityeq}
\end{align}
Here, \(\langle \cdot,\cdot\rangle\) is the inner product defined by the metric \(g\). We should note that flows obeying the two-dimensional Euler equation on a surface embedded in \(\mathbf{R}^3\) and the restriction of the foliated three-dimensional Euler flows in \(\mathbf{R}^3\) to the surface are in fact different \citep{sato2022vorticity}. Although the latter flows may be interesting from a physical point of view, we concentrate here on the pure two-dimensional Euler flows for the sake of mathematical clarity.

Let us consider the following two examples of curved surfaces. First, we take a two-dimensional sphere of unit radius as \(M\). We can take a coordinate system \((x^1, x^2) = (\lambda, \mu)\), where \(\lambda\) is the longitude and \(\mu\) is the sine of the latitude. In other words, if we embed \(M\) in \(\mathbf{R}^3\) (equipped with Euclidean metric) as \(M=\{(x,y,z)\in\mathbf{R}^3\mid x^2+y^2+z^2=1\}\), then 
\begin{align*}
    \lambda = {\rm arctan}\left(\frac{y}{x}\right), \quad \mu = z.
\end{align*}
Here, we take the Riemannian metric \(g\) of \(M\), which is inherited from \(\mathbf{R}^3\). The following equations can be shown by calculation.
\begin{align*}
    &\Gamma_{11}^1 = 0, \quad \Gamma_{11}^2 = \mu (1-\mu^2), \quad \Gamma_{12}
    ^1 = \Gamma_{21}^1 = - \frac{\mu}{1-\mu^2}, \\
    &\Gamma_{12}^2 = \Gamma_{21}^2 = 0, \quad \Gamma_{22}^1 = 0, \quad \Gamma_{22}^2 = \frac{\mu}{1-\mu^2}.
\end{align*}
On the sphere, we conventionally use the orthnormal frame \(\bm{e}_1 = (\sqrt{1-\mu^2})^{-1} \partial_1, \bm{e}_2 = \sqrt{1-\mu^2}\partial_2\) and denote the velocity as \(\bm{u}=u\bm{e}_1 + v\bm{e}_2\). In this notation, the Euler equation is expressed as:
\begin{align}
    \frac{\partial u}{\partial t} + \frac{u}{\sqrt{1-\mu^2}}\frac{\partial u}{\partial \lambda} + v\sqrt{1-\mu^2}\frac{\partial u}{\partial \mu} -\frac{\mu}{\sqrt{1-\mu^2}}uv &= -\frac{1}{\sqrt{1-\mu^2}}\frac{\partial p}{\partial \lambda}\label{k-euler-sphere-u}\\
    \frac{\partial v}{\partial t} + \frac{u}{\sqrt{1-\mu^2}}\frac{\partial v}{\partial \lambda} + v\sqrt{1-\mu^2}\frac{\partial v}{\partial \mu} + \frac{\mu}{\sqrt{1-\mu^2}} uu &= -\sqrt{1-\mu^2}\frac{\partial p}{\partial\mu}\label{k-euler-sphere-v}\\
    \frac{1}{\sqrt{1-\mu^2}}\frac{\partial u}{\partial \lambda} + \frac{\partial}{\partial \mu}(\sqrt{1-\mu^2} v) &= 0.
\end{align}
The vorticity is written as,
\begin{align*}
     q = \frac{1}{\sqrt{1-\mu^2}}\frac{\partial v}{\partial \lambda}-\frac{\partial}{\partial\mu}(\sqrt{1-\mu^2}u),
\end{align*}
and the vorticity equation is written as,
\begin{align*}
    \frac{\partial q}{\partial t} + \left(\frac{\partial \psi}{\partial \lambda}\frac{\partial q}{\partial \mu}-\frac{\partial \psi}{\partial \mu}\frac{\partial q}{\partial \lambda}\right) = 0,
\end{align*}
where \(\psi\) is the stream function obtained by inverting the elliptic relation:
\begin{align*}
q = \Delta \psi = \frac{1}{1-\mu^2}\frac{\partial^2\psi}{\partial\lambda^2} + \frac{\partial}{\partial\mu}\left((1-\mu^2)\frac{\partial \psi}{\partial\mu}\right).    
\end{align*}

The second example we consider is a two-dimensional torus \(M=(\mathbf{R}/2\pi\mathbf{Z})^2\), which is obtained by identifying points \((x,y)\) of \(\mathbf{R}^2\) with \((x+2m\pi, y+2n\pi)\in \mathbf{R}\,(m,n\in \mathbf{Z})\). We take \((x^1, x^2)=(x,y)\) as the coordinates of \(M\), and endow \(M\) with the diagonal metric:
\begin{align}
    g = \left(
    \begin{array}{cc}
        g(x,y) & 0 \\
        0 & g(x,y)
    \end{array}\right).\label{diag_riemann_metric}
\end{align}
Here, \(g(x,y)\) is a smooth function which is \(2\pi\)-periodic in both \(x\) and \(y\) directions. The coefficient of the Levi-Civita connection is given by:
\begin{align*}
        &\Gamma_{11}^1 = \Gamma_{12}^2 = \Gamma_{21}^2 =\frac{1}{2}\frac{\partial }{\partial x}\log g,\quad 
    \Gamma_{22}^2 = \Gamma_{12}^1 = \Gamma_{21}^1 = \frac{1}{2}\frac{\partial}{\partial y}\log g \\
    &\Gamma_{11}^2 = -\frac{1}{2}\frac{\partial}{\partial y}\log g,\quad 
    \Gamma_{22}^1 = -\frac{1}{2}\frac{\partial}{\partial x}\log g.
\end{align*}
As we did in the sphere example, we denote the velocity vector by \(\bm{u}=u\bm{e}_1 + v\bm{e}_2\) where \(\bm{e}_1 = (\sqrt{g})^{-1} \partial_1\) and \(\bm{e}_2 = (\sqrt{g})^{-1}\partial_2\) are the orthonormal vectors. Then the Euler equation on \(M\) is given by:
\begin{align}
    \frac{\partial u}{\partial t} + \frac{u}{\sqrt{g}}\frac{\partial u}{\partial x} + \frac{v}{\sqrt{g}}\frac{\partial u}{\partial y} + \frac{1}{2g\sqrt{g}}\frac{\partial g}{\partial y} uv - \frac{1}{2g\sqrt{g}}\frac{\partial g}{\partial x} vv &= -\frac{1}{\sqrt{g}}\frac{\partial p}{\partial x},\label{torus_eulereq_u} \\
    \frac{\partial v}{\partial t} + \frac{u}{\sqrt{g}}\frac{\partial v}{\partial x} + \frac{v}{\sqrt{g}}\frac{\partial v}{\partial y} - \frac{1}{2g\sqrt{g}}\frac{\partial g}{\partial y} uu + \frac{1}{2g\sqrt{g}}\frac{\partial g}{\partial x} uv &= -\frac{1}{\sqrt{g}}\frac{\partial p}{\partial y}, \label{torus_eulereq_v}\\
    \frac{1}{g} \left[\frac{\partial}{\partial x}(\sqrt{g}u)+\frac{\partial}{\partial y}(\sqrt{g}v)\right] &= 0. \label{torus_non_divergent}
\end{align}
The vorticity equation to be used in the following section is written as:
\begin{align}
    \frac{\partial q}{\partial t} + \frac{1}{g}\left(\frac{\partial \psi}{\partial x}\frac{\partial q}{\partial y}-\frac{\partial \psi}{\partial y}\frac{\partial q}{\partial x}\right) = 0,\label{torus_vorticity_eq}
\end{align}
where the vorticity \(q\) is defined by
\begin{align*}
    q = \frac{1}{g}\left[\frac{\partial}{\partial x}(\sqrt{g}v)-\frac{\partial}{\partial y}(\sqrt{g}u)\right],
\end{align*}
and the stream function \(\psi\) is obtained by inverting
\begin{align}
    q = \Delta \psi = \frac{1}{g}\left(\frac{\partial^2\psi}{\partial x^2} + \frac{\partial^2\psi}{\partial y^2}\right).\label{torus_qlappsi}
\end{align}
    
\subsection{Review of Haller's hyperbolic domain}
We review here the definition of hyperbolic domains introduced by \cite{haller2001lagrangian}, which was invented for two-dimensional flows and later extended to three-dimensional flows in \cite{haller2005objective}. We focus here on the case of the planar flow domain \(\mathbf{R}^2\).

Consider a material line which passes through a point \(x_0\in\mathbf{R}^2\) at time \(t=0\), and denote it by \(s\mapsto \varphi(s;0)\in\mathbf{R}^2\). Here, we take the parameter \(s\in\mathbf{R}\) so that \(\varphi(0;0)=x_0\). Let \(\varphi(s;t)\) be the position of the particle initially at \(\varphi(s;0)\) after being advected by the flow field. Furthermore, we denote the trajectory of a particle initially at \(x_0\) by \(x(t) = \varphi(0;t)\). Then, a vector field \(\bm{\xi}(t)\) along the trajectory \(x(t)\) is defined by
\begin{align*}
    \bm{\xi}(t) = \left.\frac{d}{ds}\right|_{s=0}\varphi (s; t).
\end{align*}
Note that the vector field \(\bm{\xi}(t)\) is a vector-valued function of \(t\), which can be seen as a vector field defined only on the curve \(x(t)\). 
The time evolution of \(\bm{\xi}(t)\) can be expressed by the following differential equation:
\begin{align}
    \frac{d\bm{\xi}}{dt} = (\nabla \bm{u}(x(t),t))\bm{\xi},\label{xi_evolution_R2}
\end{align}
where \(\nabla\bm{u}\) is the deformation rate tensor, defined as follows:
\begin{align*}
    \nabla \bm{u} = \left(
    \begin{array}{cc}
        \frac{\partial u}{\partial x} & \frac{\partial u}{\partial y} \\[0.25em]
        \frac{\partial v}{\partial x} & \frac{\partial v}{\partial y}
    \end{array}\right).
\end{align*}
By taking inner product of \eqref{xi_evolution_R2} and \(\bm{\xi}\), we get
\begin{align}
    \frac{1}{2}\frac{d}{dt}|\bm{\xi}|^2 = \left\langle \mathsf{S}\bm{\xi}, \bm{\xi}\right\rangle \label{xinorm_evolution_R2},
\end{align}
where \(\mathsf{S}\) is the \textit{rate-of-strain tensor}, i.e., the symmetric part of \(\nabla\bm{u}\). Differentiating \eqref{xinorm_evolution_R2} by \(t\) again, we get
\begin{align}
    &\frac{1}{2}\frac{d^2}{dt^2}|\bm{\xi}|^2 = \langle \mathsf{M}\bm{\xi},\bm{\xi}\rangle,\label{xinorm_2nddifferential_R2}\\
    &\mathsf{M} = \frac{\partial \mathsf{S}}{\partial t} + \mathsf{S}(\nabla \mathbf{u}) + {}^{t}(\nabla \mathbf{u})\mathsf{S}. \label{Mtensor_R2}
\end{align}
The tensor \(\mathsf{M}\) is called \textit{strain acceleration tensor}. If the flow field is subject to the Euler equation, we can delete the time derivative from \eqref{Mtensor_R2}, which yields:
\begin{align}
    \mathsf{M} &= \mathsf{S}(\nabla \bm{u}) + {}^{t}(\nabla \bm{u}) \mathsf{S} - \mathsf{S}^2 - \mathsf{\Omega}^2 - \mathsf{H}(p)\nonumber\\
    &= {}^{t}(\nabla \bm{u})(\nabla\bm{u}) - \mathsf{H}(p), \label{Mtensor_R2_euler}
\end{align}
where \(\mathsf{\Omega}\) is the antisymmetric part of \(\nabla\bm{u}\) and \(\mathsf{H}(p)\) is the Hesse matrix of the pressure field \(p\).

By using \eqref{xinorm_2nddifferential_R2} and \eqref{Mtensor_R2_euler}, we can diagnose the second time derivative of \(|\bm{\xi}|^2\), i.e., the acceleration of the elongation of the material line. Note that if the flow is incompressible, then \({\rm tr}\,\mathsf{S}=0\) and we can write the eigenvalues of \(\mathsf{S}\) as \(\lambda\geq 0\) and \(-\lambda\). Since the tensor \(\mathsf{S}\) is symmetric, the corresponding eigenvectors are orthogonal. In the following, we consider the generic case where \({\rm det}\,\mathsf{S}\neq 0\) at \(t=0\) (or equivalently, \(\lambda>0\) at \(t=0\)). Let \(\bm{e}_+\) be a unit eigenvector corresponding to the positive eigenvalue \(\lambda\) and let \(\bm{e}_{-}\) be a unit eigenvector corresponding to \(-\lambda\). We denote the vector \(\bm{\xi}(0)\) by using these orthonormal vectors as \(\bm{\xi}(0) = \xi^+\bm{e}_+ + \xi^-\bm{e}_-\). The equation \eqref{xinorm_evolution_R2} says that if \((\xi^+)^2-(\xi^-)^2 >0\) then \(d|\bm{\xi}|^2/dt >0\) and if \((\xi^+)^2-(\xi^-)^2 <0\) then \(d|\bm{\xi}|^2/dt<0\). Therefore, if the direction vector of the material line \(\bm{\xi}\) satisfies the first condition, the material line is elongated, and if the second condition is satisfied, the material line is contracted. The boundary of the two regimes \(Z = \{\xi^+\bm{e}_+ + \xi^-\bm{e}_- \in \mathbf{R}^2 \mid(\xi^+)^2-(\xi^-)^2 = 0\}\) is the union of two lines of \((\xi^+, \xi^-)\) plane, directed in \(\bm{\xi} = \bm{e}_+ + \bm{e}_-\) and \(\bm{\xi} = \bm{e}_+ - \bm{e}_-\). Suppose that \(\langle \mathsf{M}\bm{\xi}, \bm{\xi}\rangle \geq 0\) for two directions of \(\bm{\xi}=\bm{\xi}(0)\) in \(Z\) at \(t=0\), then we can see that there is a elongation of the fluid in the neighborhood of \(x_0\), since the slightly contracting material line will soon begin to stretch. \cite{haller2001lagrangian} demonstrated the existence of a Lagrangian coherent structure in the above condition by using the theory of stable and unstable manifolds. The \textit{hyperbolic domain} (at \(t=0\)) was defined as the set of \(x_0\in\mathbf{R}^2\) satisfying the above conditions. The \textit{elliptic domain} was defined as the complement of the hyperbolic domain. 

Note that, in three-dimensional fluid, the boundary \(Z\) of the elongating and contracting regimes is a cone. In this case, the hyperbolic and elliptic domains can be defined in the same way. However, determining them requires some effort because the cone contains infinite directions of \(\bm{\xi}(0)\) \citep{haller2005objective}. 

\subsection{Elongation of the material lines and hyperbolic domain on general two-dimensional Riemannian manifolds}

In this subsection, we derive a formula for the second time derivative of the square length of a material line element in Euler flows, in the case where the flow domain is a general orientable two-dimensional Riemannian manifold \(M\). With this formula, Haller's criterion of the hyperbolicity of the flow field can be extended to Euler flows on Riemannian manifolds.

\subsubsection{Elongation of material lines}\label{subsubsec_stretch}
As in the case of the two-dimensional plane \(\mathbf{R}^2\), which we saw in the last subsection, we can consider a material line embedded in a general two-dimensional Riemannian manifold \(M\) and the time evolution of its tangent vector \(\bm{\xi}\). Parallel to the last subsection, let \(s\mapsto \varphi(s;t)\in M\) be the material line on \(M\) and we define \(\varphi(0;0)=x_0\). The vector field \(\bm{\xi}(t)\) along the trajectory \(x(t):=\varphi(0;t)\) is defined as:
\begin{align*}
    \bm{\xi}(t) = \left.\frac{d}{ds}\right|_{s=0}\varphi(s;t).
\end{align*}
As in the case of \(\mathbf{R}^2\), \(\bm{\xi}(t)\) is defined as a vector-valued function of \(t\), and it can be seen as a vector field defined only on the curve \(x(t)\). Note that for different times \(t=t_1, t_2\, (t_1\neq t_2)\), the vectors \(\bm{\xi}(t_1)\) and \(\bm{\xi}(t_2)\) do not belong to the same tangent plane, since the flow domain \(M\) may not be flat. It may seem that the equation \eqref{xi_evolution_R2} for \(\mathbf{R}^2\) cannot be easily extended. However, we can show that it is sufficient to replace the time differentiation \(d/dt\) in \eqref{xi_evolution_R2} by the covariant derivative along the curve \(x(t)\). Then the evolution of \(\bm{\xi}(t)\) is expressed as:
\begin{align}
    \frac{D}{dt}\bm{\xi} = (\nabla_{\bm{\xi}(t)} \bm{u}) (x(t),t).\label{xi_evolution_M}
\end{align}
Here, \(D/dt\) means the covariant derivative along the curve whose parameter is \(t\), and it is not the material derivative \(D/Dt\) commonly used in fluid dynamics. The proof of \eqref{xi_evolution_M} is given in the appendix A. By taking the inner product of \eqref{xi_evolution_M} and \(\bm{\xi}\), we get
\begin{align}
        \frac{1}{2}\frac{d}{dt}|\bm{\xi}|^2 = \langle \nabla_{\bm{\xi}}\bm{u}, \bm{\xi} \rangle.\label{xinorm_evolution_M}
\end{align}

In the following, we consider an extended space \(E=M\times (-\tau, \tau)\), where \((-\tau, \tau)\subset \mathbf{R}^2\) is a space of time \(t\). For a point \((x,t)\in E\), the tangent space of \(E\) at \((x,t)\) can be written as \(T_{(x,t)}E = T_xM \oplus T_t (-\tau,\tau)\). The extended space \(E\) has a natural Riemannian metric \( g_{ij} dx^i\otimes dx^j + dt\otimes dt\). The Levi-Civita connection \(\nabla^E\) on \(E\) satisfies the following equations:
\begin{align*}
    \nabla^E_{\bm{\eta}}\bm{v} = \nabla_{\bm{\eta}}\bm{v}, \qquad \nabla^E_{\partial_t} \bm{v} = \frac{\partial \bm{v}}{\partial t}, \qquad \nabla^E_{\partial_t}\partial_t = 0,\qquad \nabla^E_{\bm{\eta}} \partial_t = 0.
\end{align*}
Here, \(\bm{\eta}\) and \(\bm{v}\) are arbitrary vector fields on \(M\) and \(\partial_t\) is a vector field corresponding to the differentiation by \(t\). The time-dependent vector field \(\bm{v}(x,t)\) on \(M\) can be viewed as a vector field on \(E\) with a vanishing component in the time direction. If we consider the curve \(t\mapsto z(t)=(x(t), t)\in E\), then \eqref{xi_evolution_M} can be written as
\begin{align}
    \nabla^E_{\dot{z}(t)}\bm{\xi} = \nabla^E_{\bm{\xi}}\bm{u}.
    \label{k1}
\end{align}
Note that, the left-hand side of (\ref{k1}) is the covariant derivative of \(\bm{\xi}\) along the curve \(z(t)\) and should be written as \(D\bm{\xi}/dt\) since it is not defined outside the curve. However, by considering an arbitrary extension \(\tilde{\bm{\xi}}\) of \(\bm{\xi}\) to the neighborhood of the curve, \(D\bm{\xi}/dt\) is equivalent to \(\nabla^E_{\dot{z}(t)}\tilde{\bm{\xi}}\), without depending on the choice of extension \(\tilde{\bm{\xi}}\) of \(\bm{\xi}(t)\). 

Before differentiating \eqref{xinorm_evolution_M} by \(t\) again, we should recall the covariant derivative of tensor fields as an element of differential geometry. First, note that the vector field \(\nabla_{\bm{\xi}}\bm{u}\) can be regarded as the value of the \((1,1)\)-tensor field \(\nabla\bm{u}\) for a vector \(\bm{\xi}\). In general, a \((1,1)\)-tensor field \(A\) defined on \(M\) means the field of linear maps \(A_x : T_x M\to T_xM\) depending smoothly on \(x\in M\). Then, the covariant derivative \(\nabla A\) is a field of linear maps \((\nabla A)_x: T_xM\otimes T_xM\to T_x M\) depending smoothly on \(x\in M\). Thus, it becomes a \((2,1)\)-tensor field, which is defined as:
\begin{align*}
    (\nabla A) (Y,X) = \nabla_Y (A(X)) - A (\nabla_Y X)
\end{align*}
for vector fields \(X\) and \(Y\) defined on \(M\). We can see that \((\nabla A)(fY, gX) = fg (\nabla A) (Y,X)\) for arbitrary smooth functions \(f,g\in C^\infty (M)\), i.e., \(\nabla A\) is a well-defined tensor field. Note that this means that the value of \((\nabla A)(Y ,X)\) at \(x\in M\) is determined only by the pointwise vectors \(X(x)\) and \(Y(x)\). Therefore, the notation \(A(Y,X)\) makes sense even if \(X\) and \(Y\) are not globally defined. It can be shown that the following formula:
\begin{align*}
    Z\langle A(X), Y\rangle = \langle (\nabla A)(Z,X), Y \rangle + \langle A (\nabla_Z X), Y \rangle + \langle A(X), \nabla_Z Y\rangle 
\end{align*}
holds for vector fields \(X, Y\) and a vector \(Z\), where the left-hand side is the derivative of \(\langle A(X), Y\rangle\) in the direction of \(Z\). 
Using this formula, we can calculate the time derivative of \eqref{xinorm_evolution_M} as:
\begin{align}
        \frac{1}{2}\frac{d^2}{dt^2}|\bm{\xi}|^2 &= \langle (\nabla^E\nabla\bm{u})(\dot{z}(t), \bm{\xi}), \bm{\xi}\rangle + \langle (\nabla \bm{u}) (\nabla^E_{\dot{z}(t)}\bm{\xi}), \bm{\xi}\rangle + \langle \nabla_{\bm{\xi}} \bm{u}, \nabla_{\dot{z}(t)}^E \bm{\xi}\rangle\nonumber\\
    & = \langle (\nabla^E\nabla\bm{u})(\bm{u}+\partial_t, \bm{\xi}), \bm{\xi}\rangle + \langle (\nabla \bm{u}) (\nabla_{\bm{\xi}}\bm{u}), \bm{\xi}\rangle + \langle \nabla_{\bm{\xi}} \bm{u}, \nabla_{\bm{\xi}} \bm{u}\rangle \nonumber\\
    & = \langle (\nabla^E\nabla\bm{u})(\partial_t, \bm{\xi}), \bm{\xi}\rangle +\langle (\nabla\nabla\bm{u})(\bm{u}, \bm{\xi}), \bm{\xi}\rangle + \langle (\nabla \bm{u}) (\nabla_{\bm{\xi}}\bm{u}), \bm{\xi}\rangle \nonumber\\
    &\qquad \qquad+ \langle \nabla_{\bm{\xi}} \bm{u}, \nabla_{\bm{\xi}} \bm{u}\rangle\nonumber\\
    &= \langle \nabla_{\bm \xi} \frac{\partial \bm{u}}{\partial t}, \bm{\xi}\rangle +\langle (\nabla\nabla\bm{u})(\bm{u}, \bm{\xi}), \bm{\xi}\rangle + \langle (\nabla \bm{u}) (\nabla_{\bm{\xi}}\bm{u}), \bm{\xi}\rangle \nonumber\\
    &\qquad \qquad+ \langle \nabla_{\bm{\xi}} \bm{u}, \nabla_{\bm{\xi}} \bm{u}\rangle.\label{xinorm_2nddiff_M_mid1}
\end{align}
To transform (\ref{xinorm_2nddiff_M_mid1}) from the third row to the fourth row, the following equation was used: 
\begin{align}
     (\nabla^E\nabla\bm{u})(\partial_t, \bm{\xi}) = \nabla_{\bm \xi} \frac{\partial \bm{u}}{\partial t}.
     \label{k2}
\end{align}
It can be shown that (\ref{k2}) holds as follows.
Since both sides of (\ref{k2}) are values of tensors, both of which have \(\bm{\xi}\) as an argument, they do not depend on the spatial distribution or time variation of \(\bm{\xi}\). Therefore, it is sufficient to show (\ref{k2}) for a fixed \((x(t_0),t_0)\in E\) and we can replace \(\bm{\xi}(t)\) by an arbitrary vector field \(\tilde{\bm{\xi}}\) defined on \(E\) which takes the value \(\bm{\xi}(t_0) \in T_{(x(t_0),t_0)}E\) at \((x(t_0),t_0)\). Consider a vector field \(\tilde{\bm{\xi}}\) that satisfies the following:
\begin{align*}
    \tilde{\bm{\xi}}(x,t) = \bar{\bm{\xi}}(x),
\end{align*}
where \(\bar{\bm{\xi}}\) is an arbitrary \(t\)-independent vector field of \(M\) which takes the value \(\bm{\xi}(t_0)\) at \(x(t_0)\in M\). Then we get the following:
\begin{align*}
     (\nabla^E\nabla\bm{u})(\partial_t, \tilde{\bm{\xi}}) &= \nabla^E_{\partial_t}(\nabla_{\overline{\bm{\xi}}}\bm{u})- \nabla_{\nabla^E_{\partial t}\overline{\bm{\xi}}}\bm{u} \\
    &= \nabla^E_{\partial_t}\nabla_{\overline{\bm{\xi}}}\bm{u} = \nabla_{\overline{\bm{\xi}}}\nabla_{\partial_t} \bm{u}= \nabla_{\bm{\xi}}\frac{\partial\bm{u}}{\partial t}.
\end{align*}
Note that we are using the fact that the vectors \(\partial_t\) and \(\bar{\bm{\xi}}\) commute.
The above shows that (\ref{k2}) holds.

Using the Euler equation \eqref{eulereq}, the time differentiation in the rightmost-hand of \eqref{xinorm_2nddiff_M_mid1} can be replaced, which gives the following:
\begin{align*}
    \frac{1}{2}\frac{d^2}{dt^2}|\bm{\xi}|^2 &= -\langle \nabla_{\bm{\xi}}({\rm grad}\,p), \bm{\xi}\rangle - \langle \nabla_{\bm{\xi}}(\nabla_{\bm{\xi}}\bm{u}), \bm{\xi}\rangle \nonumber\\
    &\qquad+ \langle (\nabla\nabla\bm{u})(\bm{u}, \bm{\xi}), \bm{\xi}\rangle + \langle (\nabla \bm{u}) (\nabla_{\bm{\xi}}\bm{u}), \bm{\xi}\rangle+ \langle \nabla_{\bm{\xi}} \bm{u}, \nabla_{\bm{\xi}} \bm{u}\rangle.
\end{align*}
Furthermore, using the identity:
\begin{align*}
    (\nabla\nabla \bm{u})(X,Y) - (\nabla\nabla \bm{u})(Y,X) = R(X,Y)\bm{u},
\end{align*}
where the \((3,1)\)-tensor field \(R: TM\otimes TM\otimes TM \to TM\) is the curvature tensor, we can derive the following simpler expression:
\begin{align}
    \frac{1}{2}\frac{d^2}{dt^2}|\bm{\xi}|^2 = -\langle \mathsf{H}(p)\bm{\xi}, \bm{\xi}\rangle - \langle R(\bm{\xi}, \bm{u})\bm{u}, \bm{\xi}\rangle + |\nabla_{\bm{\xi}} \bm{u}|^2, \label{xinorm_2nddiff_result}
\end{align}
where \(\mathsf{H}(p)\) is the Hessian tensor of the pressure \(p\), defined by \(\langle \mathsf{H}(p)\bm{\xi},\bm{\xi} \rangle = \langle \nabla_{\bm{\xi}}({\rm grad}\,p),\bm{\xi}\rangle\). Defining the strain acceleration tensor for this case as \(\mathsf{M}=-\mathsf{H}(p) -R(\cdot,\bm{u})\bm{u} + {}^{t}(\nabla\bm{u})(\nabla\bm{u})\), the right-hand side of \eqref{xinorm_2nddiff_result} is equal to \(\langle \mathsf{M}\bm{\xi},\bm{\xi}\rangle\). If \(M=\mathbf{R}^2\), we get \eqref{Mtensor_R2_euler} again.

On the right-hand side of \eqref{xinorm_2nddiff_result}, the second term reflects the curvature of \(M\). The physical meaning of this term can be clarified by using the Jacobi equation rather than fluid dynamical equations. We consider a material line \(\varphi(s;t)\) on \(M\) consisting of free mass points. Each point moves along the geodesics of \(M\) without acceleration or deceleration. The tangent vector of the material line, \(\bm{\xi}\), follows the Jacobi equation:
\begin{align*}
    \frac{D^2\bm{\xi}}{dt^2} = - R(\bm{\xi},\dot{\varphi}(s;t))\dot{\varphi}(s;t),
\end{align*}
where \(\dot{\varphi}(s;t)\) is the time derivative of \(\varphi(s;t)\). Furthermore, similar to \eqref{xi_evolution_M}, the following equation:
\begin{align*}
    \frac{D\bm{\xi}}{dt} = \nabla_{\bm{\xi}}\dot{\varphi}
\end{align*}
holds.
Therefore, using these two equations, we obtain
\begin{align*}
    \frac{1}{2}\frac{d^2}{dt^2}|\bm{\xi}|^2 &= -\langle R(\bm{\xi},\dot{\varphi})\dot{\varphi}, \bm{\xi}\rangle + \langle \nabla_{\bm{\xi}}\dot{\varphi}, \nabla_{\bm{\xi}}\dot{\varphi}\rangle,
\end{align*}
which is equivalent to \eqref{xinorm_2nddiff_result} except for the pressure term. In perfect fluids, fluid particles try to move along geodesics, but it is impossible due to the incompressibility, and the motions of the particles are corrected by pressure terms \citep{mitsumatsu2008}. In this way, we can naturally understand the meaning of \eqref{xinorm_2nddiff_result}. 

Considering the following formulas for the curvature tensor:
\begin{align*}
    &R(X,Y)Z = -R(Y,X)Z,\quad \langle R(X,Y)Z,W\rangle = -\langle R(X,Y)W, Z\rangle,\\
    &\langle R(X,Y)Z,W\rangle = \langle R(Z,W)X, Y\rangle,
\end{align*}
the curvature tensor on two-dimensional manifolds is represented by only one component: \(R_{1221} = \langle R(\bm{e}_1,\bm{e}_2)\bm{e}_2, \bm{e}_1\rangle\). Here, \(\bm{e}_1\) and \(\bm{e}_2\) are the orthonormal bases of the tangent plane of \(M\), but the value of \(R_{1221}\) does not depend on the choice of the orthonormal basis. If we set \(\bm{u} = u\bm{e}_1 + v\bm{e}_2\) and \(\bm{\xi}= \xi^1 \bm{e}_1 + \xi^2 \bm{e}_2\) in \eqref{xinorm_2nddiff_result}, then we get:
\begin{align*}
    \langle R(\bm{\xi}, \bm{u})\bm{u}, \bm{\xi}\rangle = R_{1221} (\xi^1 v - \xi^2 u)^2.
\end{align*}
Thus, at points of positive curvature (\(R_{1221}>0\)) the second term of \eqref{xinorm_2nddiff_result} is non-positive, and at points of negative curvature (\(R_{1221}<0\)), the term is non-negative. The direction of \(\bm{\xi}\) where the term is most effective is normal to \(\bm{u}\). 

\subsubsection{Extension of Haller's hyperbolicity}
    Parallel to the hyperbolic domains on a plane, we can define hyperbolic domains on a two-dimensional Riemannian manifold \(M\). For \(P\in M\), take an orthonormal basis \(\{\bm{e}_+, \bm{e}_-\}\), where \(\bm{e}_+\) and \(\bm{e}_-\) are eigenvectors of \(\mathsf{S}\) corresponding to the positive and negative eigenvalues, respectively. An arbitrary tangent vector at \(P\) can be expressed by \(\bm{\xi} = \xi^+\bm{e}_+ + \xi^- \bm{e}_-\). The hyperbolic domain on \(M\) is defined as the set of \(P\in M\) such that the strain acceleration tensor \(\mathsf{M}\) is positive definite on the set \(Z = \{\xi^+\bm{e}_+ + {\xi}^-\bm{e}_-\in T_P M\mid(\xi^+)^2-(\xi^-)^2 = 0\}\). Thus the hyperbolicity at \(P\in M\) is equivalent to \(\langle \mathsf{M}\bm{\xi}, \bm{\xi}\rangle >0\) for two vectors \(\bm{\xi} = \bm{e}_+ + \bm{e}_-\) and \(\bm{\xi} = \bm{e}_+ - \bm{e}_-\). Note that on non-flat surfaces, \({\rm tr}\,\mathsf{M}\) is not necessarily non-negative as it is on the plane \(\mathbf{R}^2\), but this is not critical for defining hyperbolicity. In fact, the existence of a Lagrangian coherent structure \citep{haller2005objective} is guaranteed by the positive definiteness of \(\mathsf{M}\) on \(Z\). The condition that \({\rm tr}\,\mathsf{M}<0\) at \(P\in M\) means that \(\mathsf{M}\) can be negative definite on \(Z\), but such a point \(P\) should be clearly excluded from the hyperbolic domain.

\section{Examples}\label{section_examples}
This section demonstrates the validity and importance of the effect expressed by \eqref{xinorm_2nddiff_result}, which implies that curvature induces acceleration or deceleration of vortex elongation, by specific examples.

\subsection{Jet on a sphere}\label{subsec_jet}
Consider the following vorticity field on a unit sphere as a simple example:
\begin{align*}
       q(\lambda,\mu) = \left\{
    \begin{array}{c}
         \frac{1}{2}(1+\mu_0)\Delta q  \qquad (\mu\geq \mu_0)\\[1.5em]
         -\frac{1}{2}(1-\mu_0)\Delta q \qquad (\mu< \mu_0),
    \end{array}
    \right.
\end{align*}
where \(\Delta q>0\) and \(\mu_0\in (-1,1)\) are constant parameters. The velocity field \(\bm{u} = u\bm{e}_1 + v\bm{e}_2\) is then given by
\begin{align*}
    u(\lambda,\mu) = \left\{
    \begin{array}{c}
         \frac{1}{2}\Delta q(1+\mu_0) \sqrt{\frac{1-\mu}{1+\mu}}\qquad (\mu\geq \mu_0)\\[1.5em]
         \frac{1}{2}\Delta q(1-\mu_0) \sqrt{\frac{1+\mu}{1-\mu}}\qquad (\mu< \mu_0).
    \end{array}
    \right. ,\qquad v(\lambda,\mu)=0.
\end{align*}
That is, the staircase vorticity field represents a jet whose axis is the latitude circle where \(\mu = \mu_0\). Considering (\ref{k-euler-sphere-v}), for this velocity field to be a steady solution, the following equation should hold:
\begin{align*}
    -\sqrt{1-\mu^2}\frac{\partial p}{\partial \mu} = \frac{\mu}{\sqrt{1-\mu^2}}u^2.
\end{align*}
Thus, for \(\mu > \mu_0\), the Hesse tensor of the pressure field \(p\) is calculated as:
\begin{align*}
        H(p) = \left(
    \begin{array}{cc}
        \frac{\mu^2}{(1+\mu)^2}(1+\mu_0)^2 \left(\frac{\Delta q}{2}\right)^2 & 0 \\[1.0em]
        0 & \frac{2\mu-1}{(1+\mu)^2}(1+\mu_0)^2 \left(\frac{\Delta q}{2}\right)^2
    \end{array}
    \right).
\end{align*}
Also in this case, since the deformation velocity tensor is calculated as follows:
\begin{align*}
    \nabla \bm{u} = \left(
    \begin{array}{cc}
        0 & -\frac{1}{1+\mu}(1+\mu_0)\frac{\Delta q}{2} \\[1em]
        \frac{\mu}{1+\mu}(1+\mu_0)\frac{\Delta q}{2} & 0
    \end{array}
    \right),
\end{align*}
we finally get
\begin{align*}
    -H(p) + {}^{\rm T}(\nabla \bm{u}) (\nabla\bm{u}) = \left(
    \begin{array}{cc}
       0  & 0 \\[1em]
       0  & \frac{2(1-\mu)}{(1+\mu)^2}(1+\mu_0)^2\left(\frac{\Delta q}{2}\right)^2
    \end{array}
    \right).
\end{align*}
Furthermore, since \(R_{1221}=1\) anywhere on the sphere, we can write as:
\begin{align*}
    R(\cdot, \bm{u})\bm{u} = \left(
    \begin{array}{cc}
       v^2  & -uv \\[1em]
       -uv  &  u^2
    \end{array}\right) =
    \left(
    \begin{array}{cc}
        0 & 0 \\[1em]
        0 & \frac{1-\mu}{1+\mu}(1+\mu_0)^2\left(\frac{\Delta q}{2}\right)^2
    \end{array}
    \right).
\end{align*}
Thus, for a tangent vector \(\bm{\xi} = \xi^1\bm{e}_1 + \xi^2\bm{e}_2\) of the material line, we have
\begin{align}
    \frac{1}{2}\frac{d^2}{dt^2}|\bm{\xi}|^2 &= -\langle R(\bm{\xi},\bm{u})\bm{u}, \bm{\xi}\rangle - \langle H(p)\bm{\xi},\bm{\xi}\rangle + \langle \nabla_{\bm{\xi}}\bm{u}, \nabla_{\bm{\xi}}\bm{u}\rangle \nonumber\\
    &= \left[-\frac{1-\mu}{1+\mu}(1+\mu_0)^2 + \frac{2(1-\mu)}{(1+\mu)^2}(1+\mu_0)^2\right]\left(\frac{\Delta q}{2}\right)^2 (\xi^2)^2\nonumber\\
    &= \frac{(1-\mu)^2}{(1+\mu)^2}(1+\mu_0)^2 \left(\frac{\Delta q}{2}\right)^2 (\xi^2)^2.\label{example1_length}
\end{align}
Note that the length of a material line can be calculated by integration because in this case the trajectory of a particle can be written explicitly. If the initial position of an element of the material line is given by \(\varphi(s; 0) = (\lambda(s; 0), \mu(s; 0)) = (0, z_0 + s\Delta \mu)\, (s\in[0,1])\), then the position at time \(t\) can be given by \(\varphi(s;t)=(\lambda (s;t), \mu (s;t))\), where each component is written as:
\begin{align*}
     \lambda(s; t) = \frac{1}{2}\Delta q(1+\mu_0) \sqrt{\frac{1-(z_0 + s\Delta \mu)}{1+(z_0 + s\Delta \mu)}} t,\qquad \mu(s;t) = z_0 + s\Delta \mu.
\end{align*}
Thus, the length \(l(t)\) of the material line at time \(t\) can be calculated as:
\begin{align*}
        l(t) = \int_{z_0}^{z_0+\Delta \mu} \sqrt{\left(\frac{1-\mu}{1+\mu}\right)^2 \left[\frac{\Delta q}{2}(1+\mu_0) t\right]^2 + \frac{1}{1-\mu^2}} d\mu.
\end{align*}
If \(|\Delta \mu|\ll 1\), we have
\begin{align*}
    l(t) \approx \sqrt{\left(\frac{1-z_0}{1+z_0}\right)^2 \left[\frac{\Delta q}{2}(1+\mu_0) t\right]^2 + \frac{1}{1-z_0^2}} \Delta\mu
\end{align*}
and 
\begin{align*}
    \frac{1}{2}\frac{d^2}{dt^2}|l(t)|^2 \approx \left(\frac{1-z_0}{1+z_0}\right)^2 (1+\mu_0)^2 \left(\frac{\Delta q}{2}\right)^2 |\Delta \mu|^2.   
\end{align*}
Therefore, \eqref{example1_length} is consistent with the explicit expression of the length of the material line.

\subsection{Hyperbolic domains of vortices on a sphere}\label{subsec_hypbsphere}
In the case of the zonal jet shown above, the strain acceleration tensor \(\mathsf{M}\) is positive semidefinite.  The entire sphere therefore belongs to the hyperbolic domain in Haller's sense, except for the latitude circle \(\mu=\mu_0\) where the rate-of-strain tensor \(\mathsf{S}\) is zero. This is the case for all zonal shear flows on the sphere, but for non-zonal flows the separation between hyperbolic and elliptic is non-trivial. In particular, the effect of the curvature cannot be neglected in determining a hyperbolic domain.

\begin{figure}[h]
    \centering
    \includegraphics[width=0.65\linewidth]{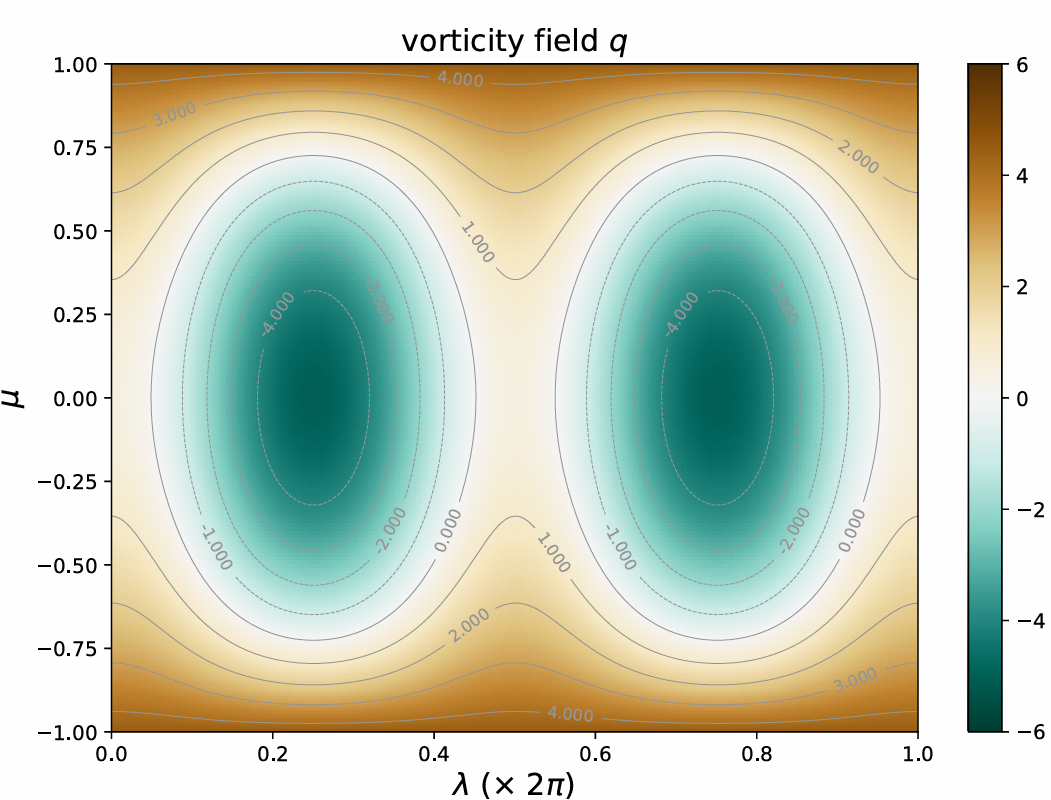}
    \vspace{6em}
    \caption{The vorticity field \(q\) defined by \eqref{sphere_q_def}. The horizontal axis is longitude \(\lambda\), and the vertical axis is sine latitude \(\mu\).}
    \label{fig:sphere_q}
\end{figure}

\begin{figure}[htbp]
    \centering
    \includegraphics[width=0.65\linewidth]{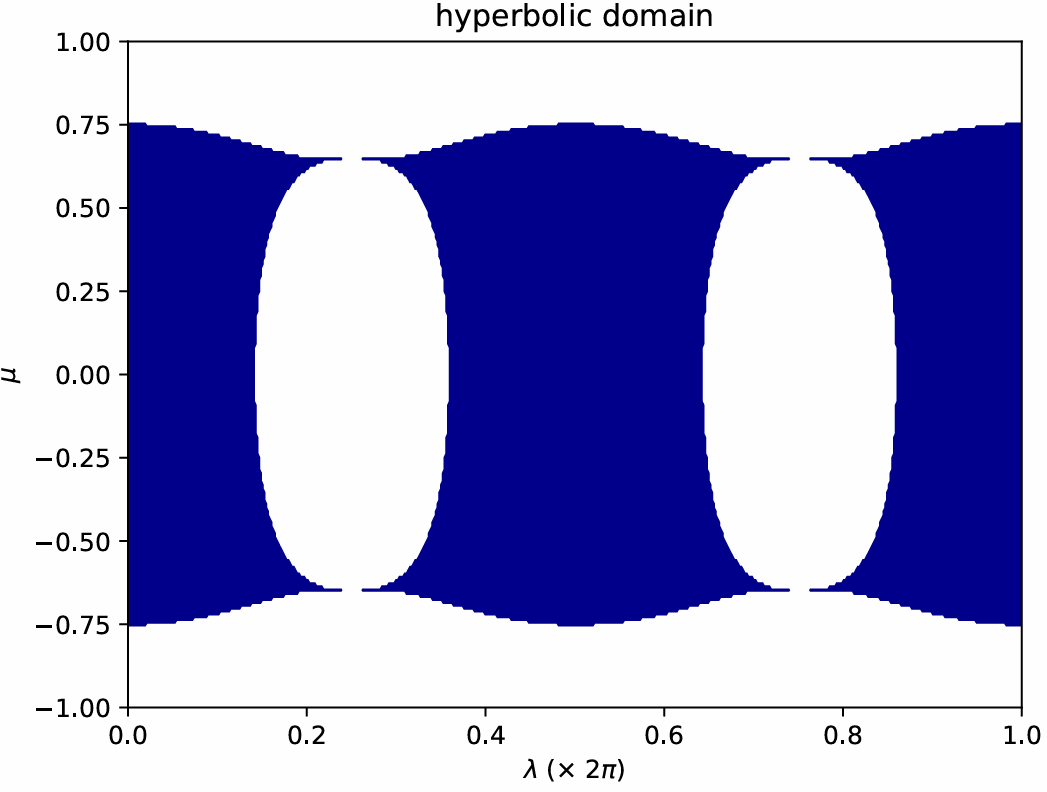}\\
    \includegraphics[width=0.65\linewidth]{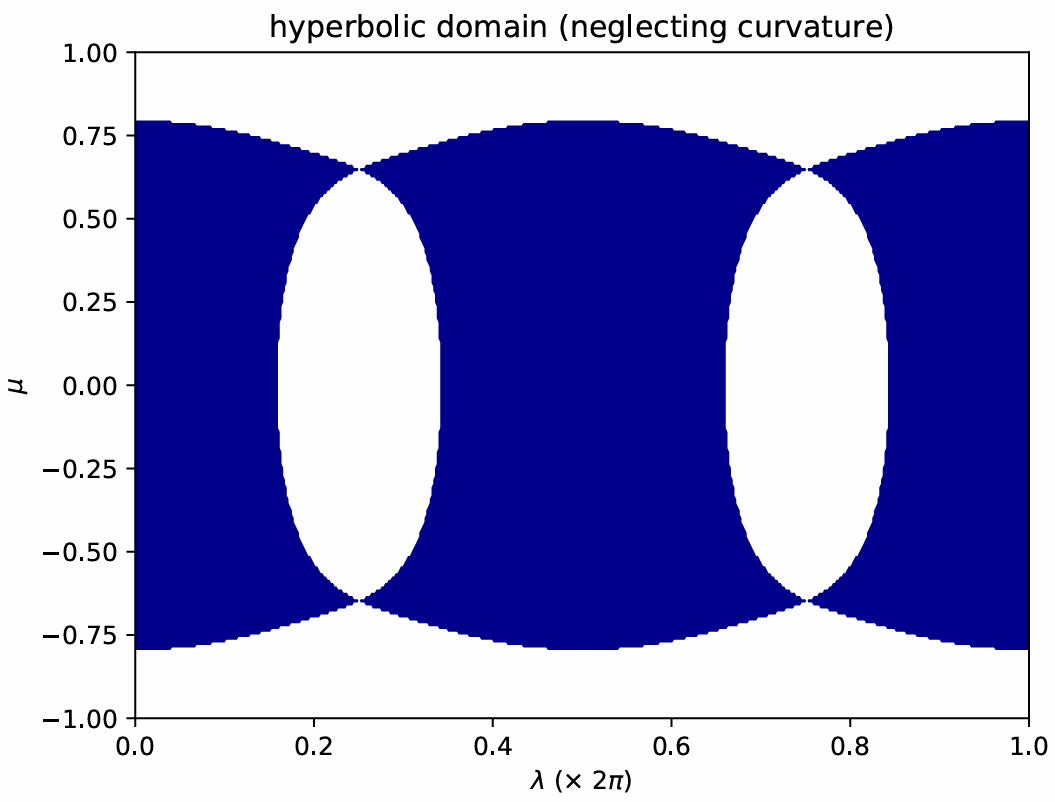}
    \vspace{12em}
    \caption{Top panel: the hyperbolic domain (colored blue) of the vorticity field \eqref{sphere_q_def}. Bottom panel: the same as the top panel, except that the curvature term of the strain acceleration tensor is neglected.}
    \label{fig:hypb_sphere}
\end{figure}

Here we consider a steady vorticity field on a unit sphere (figure \ref{fig:sphere_q}), defined as:
\begin{align}
    q(\lambda,\mu) = \frac{\sqrt{5}}{2}(3\mu^2-1) + \frac{\sqrt{30}}{2}(1-\mu^2)\cos 2\lambda.\label{sphere_q_def}
\end{align}
The first and second terms on the right-hand side of \eqref{sphere_q_def} are eigenvectors of the Laplace-Beltrami operator on the sphere with the common eigenvalue \(-6\). In fact, these terms are twice the real parts of the spherical harmonics \(Y_{0,2}(\lambda,\mu)\) and \(Y_{2,2}(\lambda,\mu)\), respectively. The numerically computed hyperbolic domain for the vorticity field \eqref{sphere_q_def} is shown in the top panel of figure \ref{fig:hypb_sphere}. On the other hand, the hyperbolic domain calculated by deliberately neglecting the curvature term of the strain acceleration tensor is shown in the bottom panel of figure \ref{fig:hypb_sphere}. Since the curvature term is negative semidefinite on the sphere, the hyperbolic domain with the curvature term properly considered is smaller than that with the curvature term neglected.

\begin{figure}
    \centering
    \includegraphics[width=0.65\linewidth]{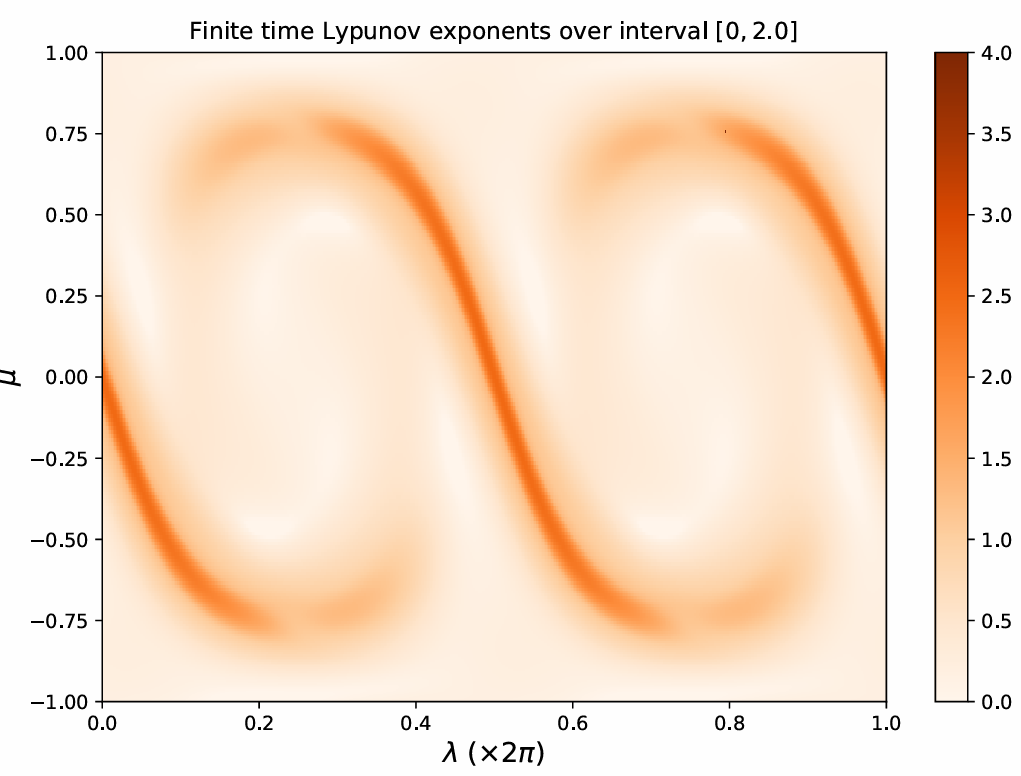}
    \vspace{6em}
    \caption{The (forward) finite-time Lyapunov index over the time interval from \(t=0\) to \(t=2\) for the flow given by the vorticity field defined by \eqref{sphere_q_def}.}
    \label{fig:hypb_FTLE_sphere}
\end{figure}

\begin{figure}
    \centering
    \includegraphics[width=0.65\linewidth]{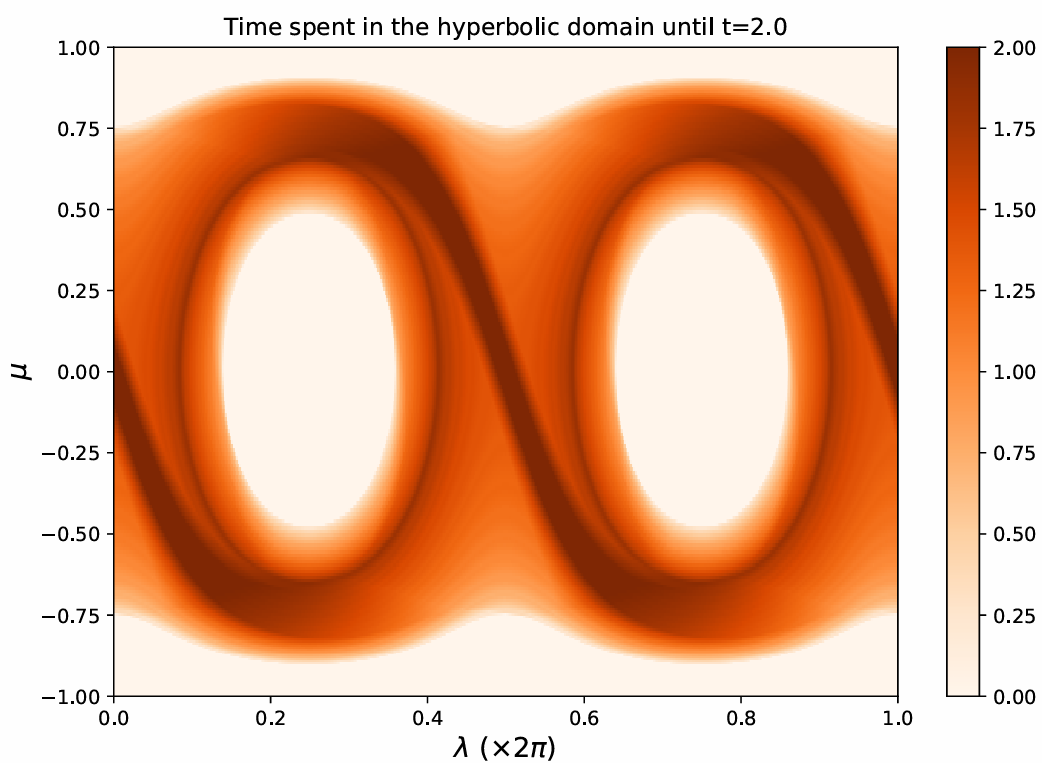}
    \includegraphics[width=0.65\linewidth]{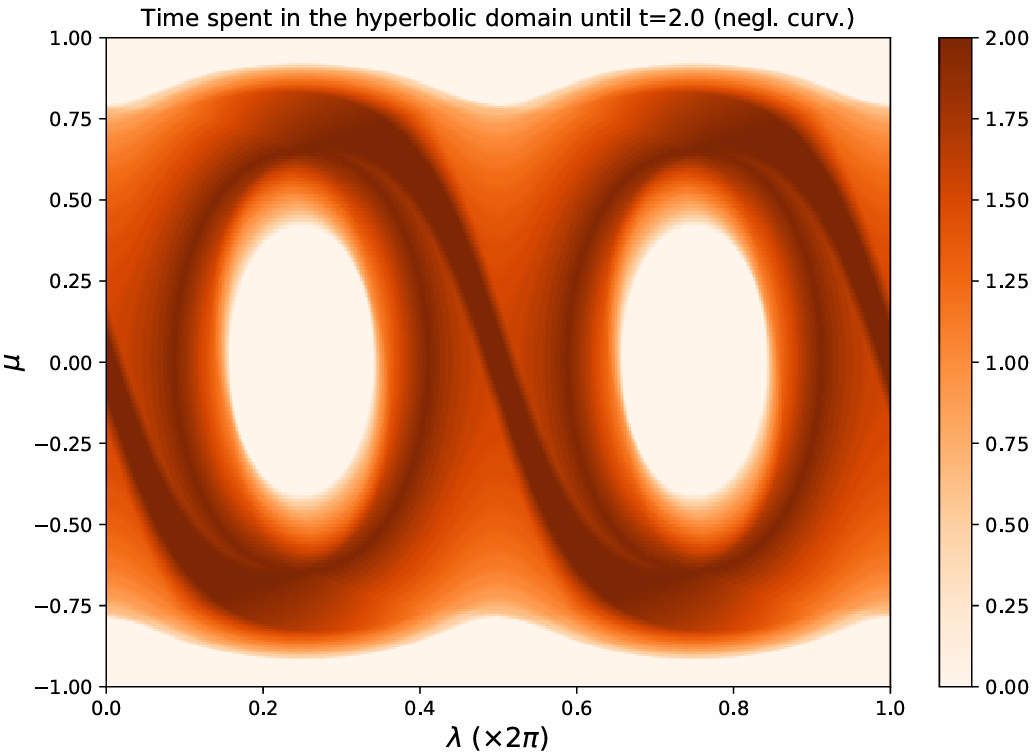}
    \vspace{6em}
    \caption{Top panel: total time spent by a particle starting from each point at \(t=0\) in the hyperbolic domain up to \(t=2\) for the flow field given by \eqref{sphere_q_def}. Bottom panel: same as the top panel, except that the hyperbolic domain is computed by neglecting the curvature term of the strain acceleration tensor \(\mathsf{M}\).}
    \label{fig:hypbtime}
\end{figure}

Following the analysis of \cite{haller2005objective} for the Arnold-Beltrami-Childress flow in a three-dimensional flat torus, we computed the finite-time Lyapunov index for the vorticity field \eqref{sphere_q_def}. The finite-time Lyapunov index for the time interval from \(t=0\) to \(t=2\) is shown in figure \ref{fig:hypb_FTLE_sphere}. The stream function for the vorticity field has two saddles at \((\lambda,\mu)=(0,0)\) and \((\pi,0)\). At each of these two saddles, a repelling material line and an attracting material line cross, both of which are time-invariant. As shown in figure \ref{fig:hypb_FTLE_sphere}, the finite-time Lyapunov index is large in the band-like regions containing the two saddles, which correspond to the repelling material lines.  The top panel of figure \ref{fig:hypbtime} shows the time that a particle leaving each point at \(t = 0\) spends in the hyperbolic domain until \(t = 2\) (hyperbolicity time). In the neighborhood of the repelling material lines with large values of the finite-time Lyapunov index, the hyperbolicity time is also large. However, there are annular regions where the hyperbolicity time has a local maximum along the negative vortices. Each of these regions corresponds to a periodic orbit of fluid particles passing the cusp of the hyperbolic domain (\((\lambda,\mu)\approx (\pi/4,0.65)\) and \((3\pi/4,0.65)\)). Although a small material line on this periodic orbit spends a long time in the hyperbolic domain, it is not particularly elongated compared to other material lines starting from a point outside the annular maximum. In fact, such a material line is gradually elongated as it repeats the oscillatory elongation and contraction on the periodic orbits, and the total elongation rate after a certain time is not significantly large compared to material lines traveling outside the periodic orbit. We should note that even in the hyperbolic domain, the right-hand side of the equation: \(d
^2|\bm{\xi}|^2/dt^2 = 2\langle \mathsf{M}\bm{\xi},\bm{\xi}\rangle\) is not necesarilly positive, depending on the direction of \(\bm{\xi}\). The annular regions of large hyperbolicity time become more pronounced when we use the hyperbolicity time with the curvature term neglected (bottom panel of figure \ref{fig:hypbtime}). 

\begin{figure}
    \centering
    \includegraphics[width=0.65\linewidth]{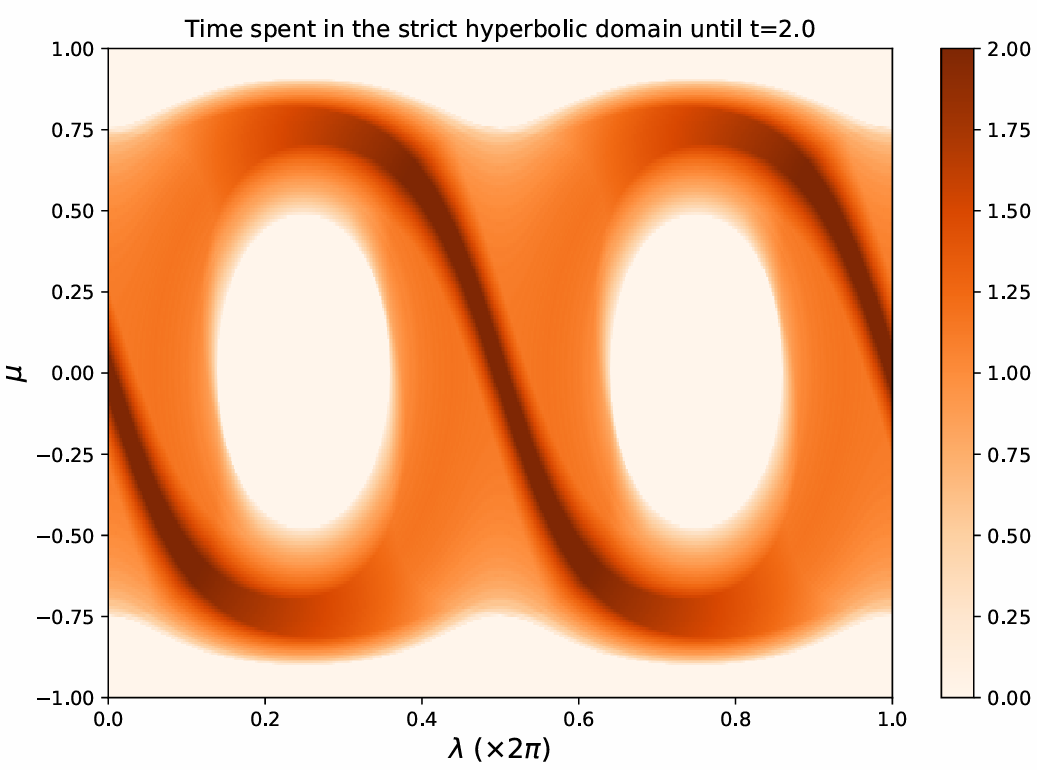}
    \vspace{6em}
    \caption{Total time spent by a particle starting from each point at \(t=0\) in the \textit{strong} hyperbolic domain up to \(t=2\) for the flow field given by \eqref{sphere_q_def}.}
    \label{fig:strict_hypbtime}
\end{figure}

For steady vorticity fields with periodic orbits like this example, it may be more appropriate to use the \textit{strong hyperbolic domain} \citep{haller2005objective} rather than the normal hyperbolic domain. The strong hyperbolic domain is the set of points where the strain acceleration tensor \(\mathsf{M}\) is positive definite. In the strong hyperbolic domain, \(d^2|\bm{\xi}|^2/dt^2\) is always positive, regardless of the direction of \(\bm{\xi}\). The time evolution of \(|\bm{\xi}(t)|^2\) is always convex as long as the material line is included in the strong hyperbolic domain and no oscillation occurs. Figure \ref{fig:strict_hypbtime} shows the total time spent by a particle starting from each point at \(t=0\) in the strong hyperbolic domain up to \(t=2\). Note that the annular maximum of hyperbolicity time does not appear.

\subsection{Onset of the vorticity filamentation on a curved torus}\label{subsec_torus}
So far we have considered the flows on a sphere, where the curvature \(R_{1221}\) is a positive constant and the curvature term of the strain acceleration tensor acts to decelerate the elongation of material lines. Conversely, if we consider a surface with regions of negative curvature, the elongation of a material line will be accelerated at points of negative curvature. In Euler flows, the elongation of a material line can lead to vorticity filamentation because the vorticity is frozen to the fluid particles. Thus, the negative curvature of the surface itself can promote the filamentation of the vorticity field. In this subsection, we perform a numerical time integration of the Euler equation on a curved torus and check the effect of negative curvature on vorticity filamentation.

\begin{figure}
    \centering
    \includegraphics[width=0.5\linewidth]{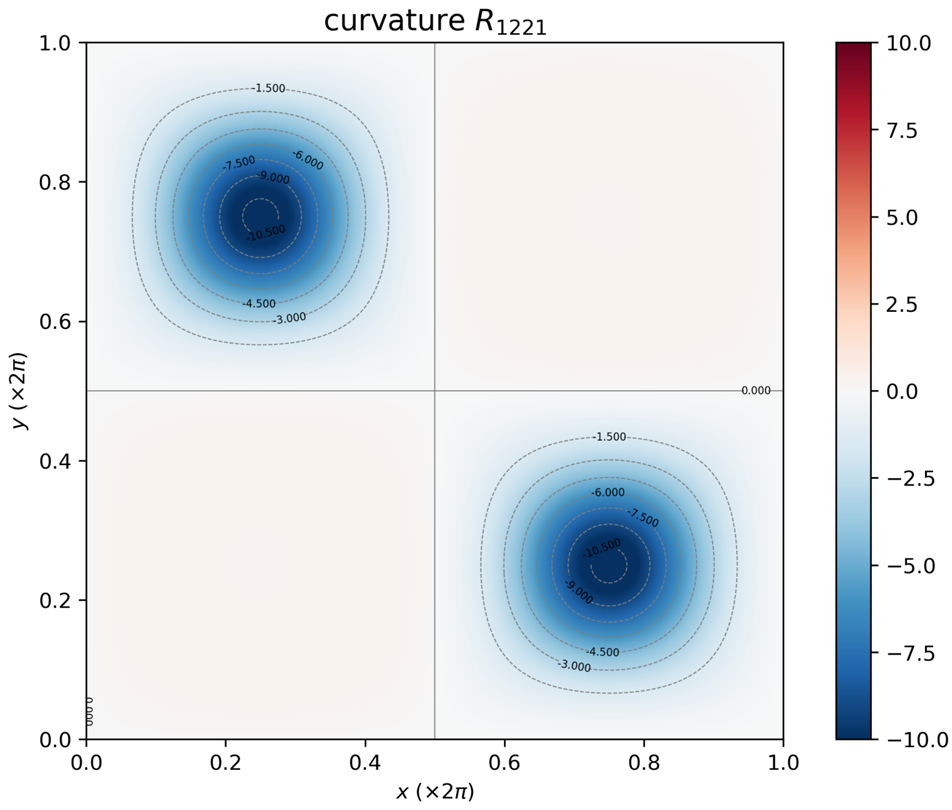}
    \caption{The value of the curvature \(R_{1221}\) on the torus with the setting \(\alpha = 1.8\). The horizontal axis is \(x\), and the vertical axis is \(y\). Note that the maximum value of \(R_{1221}\) is \(\alpha e^{-\alpha} \approx 0.29754\), and the minimum value of \(R_{1221}\) is \(-\alpha e^{\alpha}\approx -10.889\).}
    \label{fig:torus_curvature}
\end{figure}

\begin{figure}
    \centering
    \includegraphics[width=0.65\linewidth]{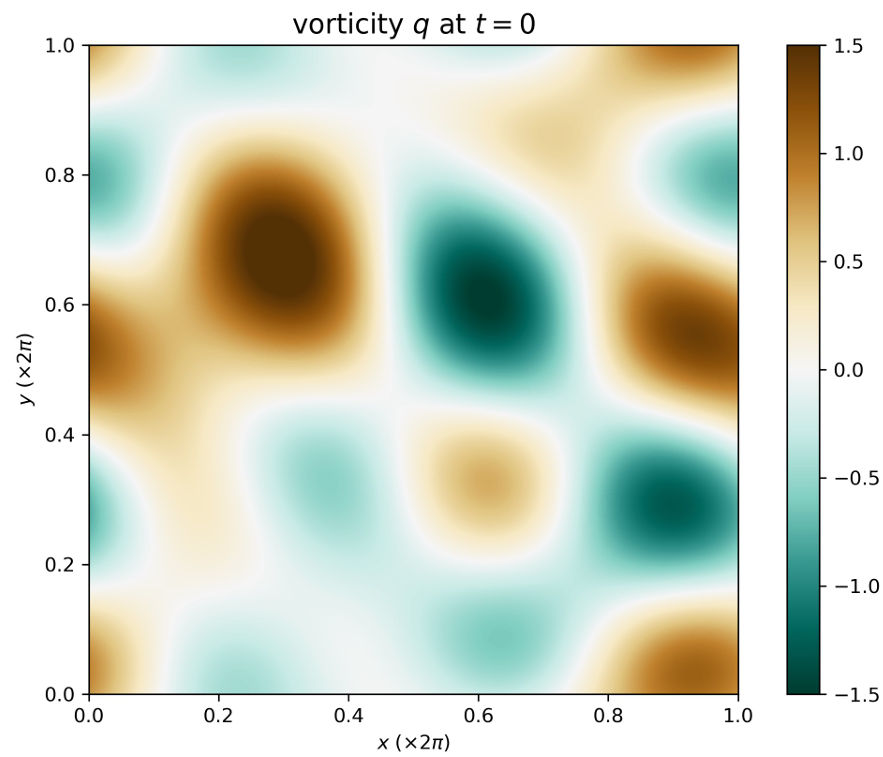}
    \caption{The initial vorticity field \(q\) defined on the curved torus. The horizontal axis is \(x\), and the vertical axis is \(y\).}
    \label{fig:torus_initvol}
\end{figure}

\begin{figure}
    \centering
    \begin{tabular}{cc}
        \begin{minipage}[t]{0.49\hsize}
            \includegraphics[scale=0.30]{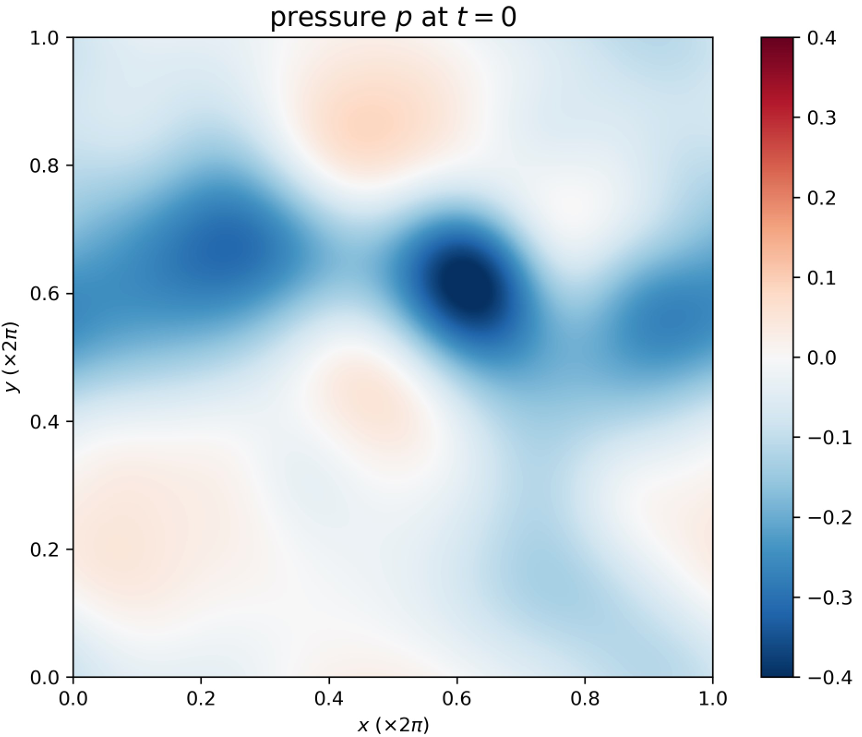}       
            \centering
        \end{minipage} 
        &
        \begin{minipage}[t]{0.49\hsize}
            \includegraphics[scale=0.30]{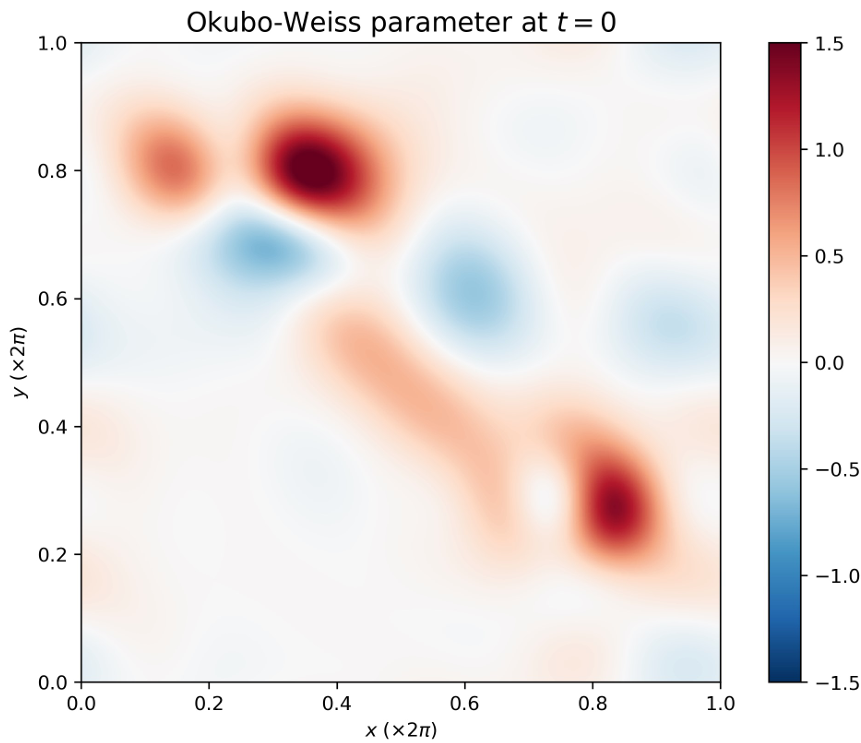}    
            \centering
        \end{minipage}
    \end{tabular}
    \caption{Left panel: the pressure field \(p\) for the initial vorticity field shown in figure \ref{fig:torus_initvol}. Right panel: The Okubo-Weiss parameter \(Q\) for the same initial vorticity field.}
    \label{fig:torus_init_press_okuboweiss}
\end{figure}

We consider a two-dimensional torus, which has a Riemannian metric defined by \eqref{diag_riemann_metric} with the following explicit form of function:
\begin{align*}
    g(x,y) = e^{\alpha\sin x\sin y},
\end{align*}
where \(\alpha >0\) is a constant. The curvature of the torus is then calculated as:
\begin{align*}
    R_{1221} &= -\frac{1}{2}g^{-1} \left(\frac{\partial^2}{\partial x^2}+ \frac{\partial^2}{\partial y^2}\right)\log g \\
    &= \alpha \sin x \sin y \cdot e^{-\alpha \sin x\sin y}.
\end{align*}
Therefore, the curvature is positive at points where \(\sin x \sin y >0\) and negative at points where \(\sin x \sin y <0\). We set \(\alpha =1.8\) for this example, and then the distribution of the curvature \(R_{1221}\) is shown in figure \ref{fig:torus_curvature}. The initial vorticity field \(q\) is given by a linear combination of the trigonometric functions \(\exp(\sqrt{-1}(kx + ly))\) of \(|k|\leq 2\) and \(|l|\leq 2\) with randomly generated coefficients. Note that the constant part of the initial vorticity is adjusted so that the area mean (weighted by \(g(x,y)\)) of \(q\) is zero. The explicit definition of the initial vorticity field is given in the appendix C for the sake of reproducibility. The initial vorticity field is shown in figure \ref{fig:torus_initvol}. Note that the area element on the torus is \(g dx dy\), and therefore, regions with large \(g\) appear to be shrunk, and regions with small \(g\) appear to be expanded in the figure. The corresponding pressure field \(p\) is shown in the left panel of figure \ref{fig:torus_init_press_okuboweiss}. The right panel of figure \ref{fig:torus_init_press_okuboweiss} shows the distribution of the Okubo-Weiss parameter, defined as:
\begin{align*}
    Q = \frac{1}{2}(|\mathsf{S}|^2-|\mathsf{\Omega}|^2).
\end{align*}
As shown in figure \ref{fig:torus_initvol}
there is a positive vorticity region near \((x,y)=(0.30\times 2\pi, 0.70\times 2\pi)\). The pressure around this region is negative. Furthermore, the Okubo-Weiss parameter is negative in the core of this region as shown in figure \ref{fig:torus_init_press_okuboweiss}, which means that the positive vorticity region is an Eulerian elliptic domain and can be considered as a vortex.

\begin{figure}
    \centering
    \includegraphics[width=1.0\linewidth]{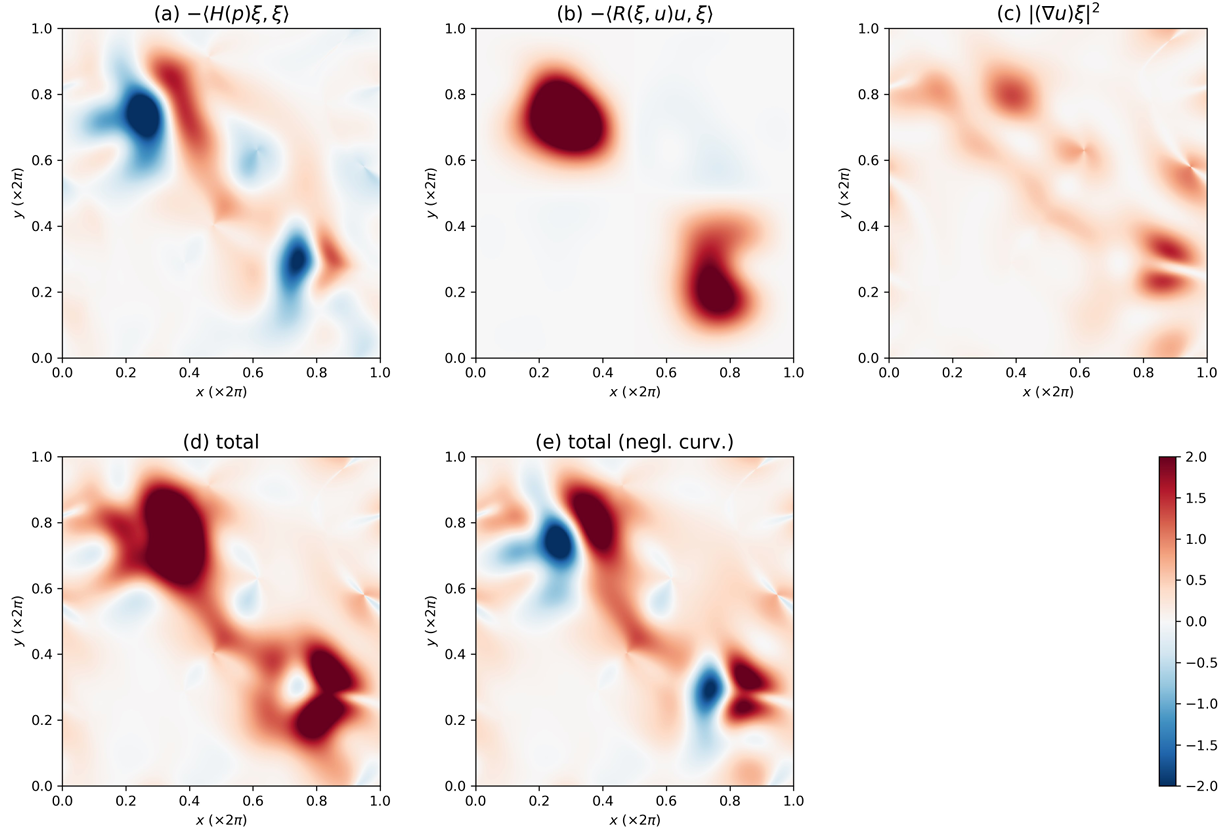}
    \caption{For the initial vorticity field shown in figure \ref{fig:torus_initvol}, the terms of the strain acceleration \(\langle \mathsf{M}\bm{\xi}, \bm{\xi}\rangle\) and their sum with setting \(\bm{\xi}\) to be the unit vector normal to the velocity vector \(\bm{u}\) at each point. (a) \(-\langle \mathsf{H}(p)\bm{\xi},\bm{\xi}\rangle\): the Hesse tensor term of the pressure \(p\), (b) \(-\langle R(\bm{\xi}, \bm{u})\bm{u},\bm{\xi}\rangle\): the curvature term, (c) \(\langle {}^{t}(\nabla \bm{u})(\nabla\bm{u})\bm{\xi}, \bm{\xi}\rangle = |\nabla_{\bm{\xi}}\bm{u}|^2\): the term due to the velocity gradient, (d) \(\langle \mathsf{M}\bm{\xi},\bm{\xi}\rangle\): the sum, and (e) \(\langle \mathsf{M}\bm{\xi}, \bm{\xi}\rangle - \langle R(\bm{\xi},\bm{u})\bm{u},\bm{\xi}\rangle\): the sum without the curvature term.}
    \label{fig:torus_hypb_term}
\end{figure}

The value of each term on the right-hand side of \eqref{xinorm_2nddiff_result} computed for the initial vorticity field is shown in figure \ref{fig:torus_hypb_term}. Here we choose a unit vector normal to the pointwise velocity as \(\bm{\xi}\) at each point. The total value \(\langle \mathsf{M}\bm{\xi},\bm{\xi}\rangle\) and the total value without the curvature term are also shown. Note that in the negative curvature region, there are points where \(d^2|\bm{\xi}|^2/dt^2\) is positive, although it would be negative if we neglected the effect of curvature. 

The time integration of the vorticity equation was performed by the spectral method using the usual trigonometric function expansion of \(x\) and \(y\). Note that in the computation from/to the stream function \(\psi\) to/from the vorticity \(q\), or in the computation of the nonlinear term of the vorticity equation, we have to compute multiplications involving the Riemannian metric \(g(x,y)\) or its inverse \(\{g(x,y)\}^{-1}\). The computation of these multiplications was done in the physical space, and the results were transformed into the spectral space. Although \(g\) and \(g^{-1}\) have infinite Fourier series expansions, we do not have to worry much about aliasing errors because the high wavenumber components decay rapidly. In fact, for sufficiently large \(k\) and \(l\), the spectral coefficients \(a_{k,l}\) of \(g\) or \(g^{-1}\) are evaluated as follows:
\begin{align}
    |a_{k,l}| \leq C\frac{1}{L(k,l)!}\left(\frac{\alpha}{2}\right)^{L(k,l)},\label{metric_g_evaluation}
\end{align}
where \(C>0\) is a constant that does not depend on \(k\) and \(l\), and \(L(k,l) := \max \{|k|,|l|\}\). The proof of this inequality is given in the appendix B. The decay is particularly fast for \(\alpha<2\), so we chose \(\alpha = 1.8\) instead of \(\alpha =2\). The time integration was done using \(800\times 800\) grids in physical space with equal spacing. The truncation wavenumber was set to \(K=199\). We added a small viscosity term on the right-hand side of the vorticity equation to avoid the artificial accumulation of enstrophy at the grid scale. That is, we integrated the following vorticity equation:
\begin{align*}
    \frac{\partial q}{\partial t} + \frac{1}{g}\left(\frac{\partial \psi}{\partial x}\frac{\partial q}{\partial y}- \frac{\partial \psi}{\partial y}\frac{\partial q}{\partial x}\right) = \nu \frac{1}{g}\left(\frac{\partial^2 q}{\partial x^2}+\frac{\partial^2 q}{\partial y^2}\right).
\end{align*}
Note that we should use the form of the viscosity term that takes into account the Ricci curvature of the manifold, which is consistent with the Navier-Stokes equation on Riemannian manifolds, as done in \cite{reuther2015interplay}. However, since we do not continue the time integration until the viscosity term plays a crucial role, the choice of specific forms of the viscosity term is not important here, and we therefore use the simplest form of the viscosity term to remove the accumulation of enstrophy at high wavenumbers. The value of the viscosity coefficient \(\nu\) was set to \(\nu=1.0\times 10^{-6}\). Time marching was performed using the fourth-order Runge-Kutta method with the time step of \(\Delta t = 0.001\).

The result of the time integration is shown in figure \ref{fig:torus_evolution} (the gif animation file is available in the supplementary data). The positive vortex region, initially located around the point \((x,y)=(0.30\times 2\pi, 0.70\times 2\pi)\), is deformed over time. The positive vortex is elongated and bent up to the time \(t = 3.0\). The vorticity of the bent region is then concentrated around \((x, y) = (0.20 \times 2\pi,0)\), and the rest of the vortex is deformed into a filamentary shape. Note that only \(1.87\times 10^{-3} \)\% of the initial energy and \(4.03\times 10^{-3}\)\% of the initial enstrophy are lost through dissipation by the time \(t=5\). 

\begin{figure}
    \centering
    \includegraphics[width=1.0\linewidth]{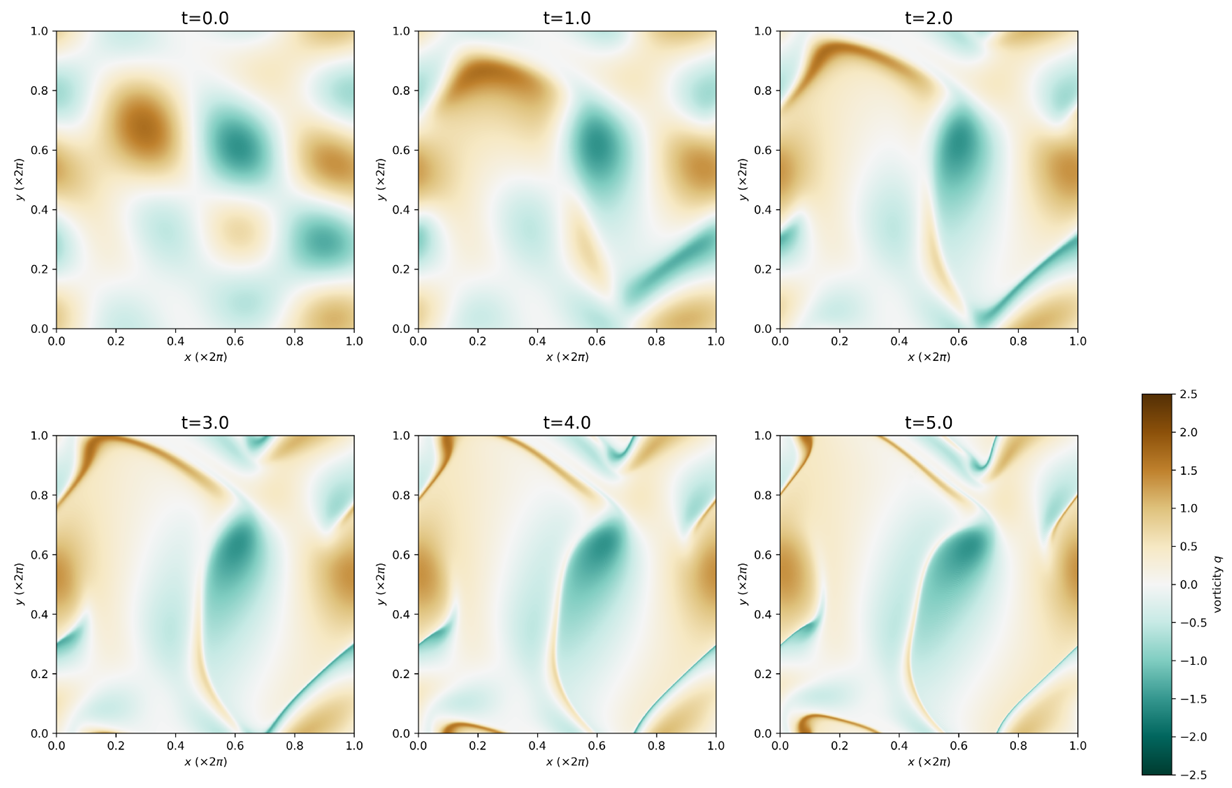}
    \caption{The time evolution of the vorticity field \(q\) from the initial vorticity field shown in figure \ref{fig:torus_initvol}.}
    \label{fig:torus_evolution}
\end{figure}

\begin{figure}
    \centering
    \includegraphics[width=1.0\linewidth]{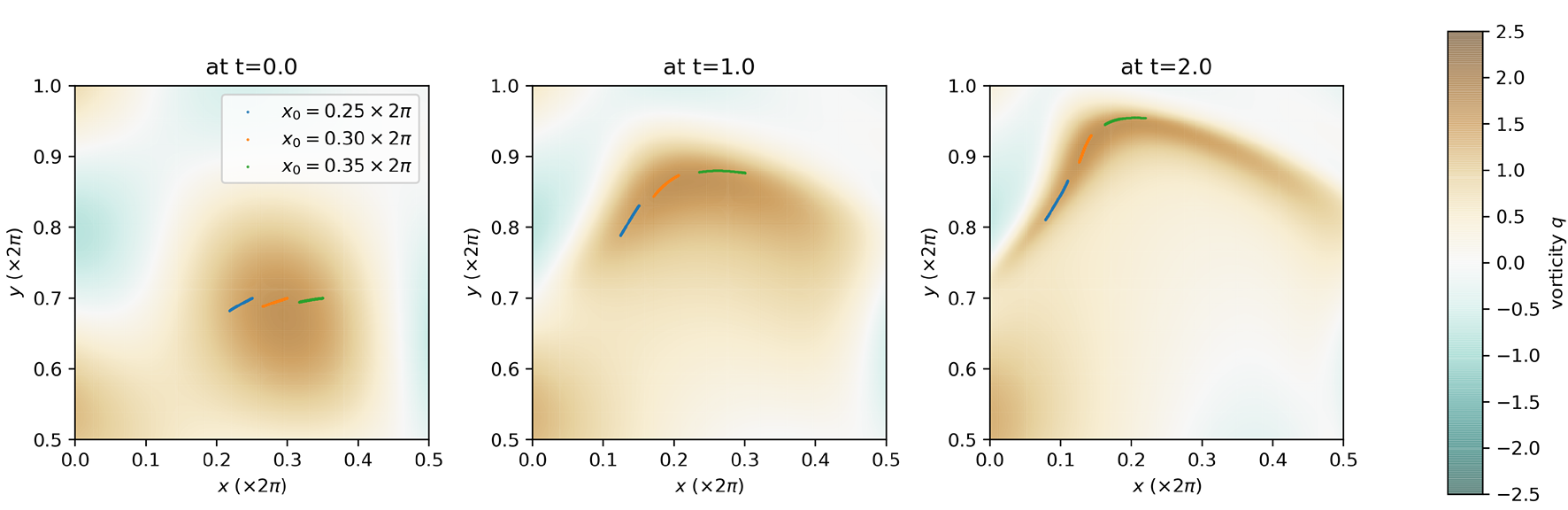}
    \caption{The position of the nodes of the material lines at \(t=0, 1\), and \(t=2\). The nodes are drawn in blue for \((x_0, y_0) =(0.25\times 2\pi, 0.70\times 2\pi)\), orange for \((x_0, y_0) =(0.30\times 2\pi, 0.70\times 2\pi)\), and green for \((x_0, y_0) =(0.35\times 2\pi, 0.70\times 2\pi)\). The time \(t\) is shown at the top of each panel, and the vorticity field at \(t\) is also shown. Note that only the part of the torus \(0\leq x\leq 0.5\times 2\pi, 0.5\times 2\pi\leq y\leq 2\pi\) is shown.
    }
    \label{fig:torus_materialline}
\end{figure}

\begin{figure}
    \centering
    \includegraphics[width=0.6\linewidth]{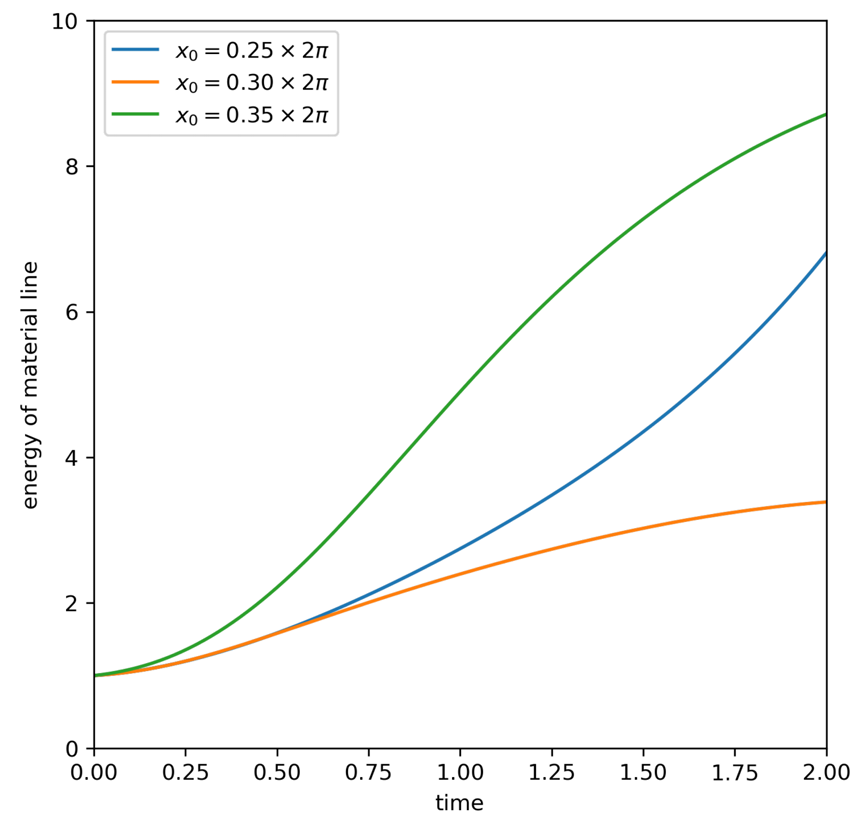}
    \caption{The time evolution of the energy \(\mathcal{E}\) of each material line in \(t\in[0,2]\). The horizontal axis is the time \(t\). The vertical axis is the energy \(\mathcal{E}(t)\) of the material line normalized by the initial energy \(\mathcal{E}(0)\).
    }
    \label{fig:torus_materialline_energy}
\end{figure}

To quantitatively confirm the elongation of the vortex, we compute the advection of material lines. We defined an initial material line \(\varphi(s;0) = (x(s), y(s))\) as the solution of the following differential equation:
\begin{align*}
    &\frac{dx}{ds} = -\frac{1}{\sqrt{g}}\frac{v}{\sqrt{u^2+v^2}}, \quad \frac{dy}{ds} = \frac{1}{\sqrt{g}}\frac{u}{\sqrt{u^2+v^2}},\\
    & x(0) = x_0, y(0) = y_0, s\in [0, 0.1].
\end{align*}
Here, \(u\) and \(v\) are the components of the velocity vector \(\bm{u} = u\bm{e}_1 + v\bm{e}_2\) at \(t=0\). By definition, the material line is normal to the velocity vectors at each point. For the three initial points: \((x_0,y_0) = (0.25\times 2\pi, 0.70\times 2\pi), (0.30\times 2\pi, 0.70\times 2\pi)\) and \((0.35\times 2\pi, 0.70\times 2\pi)\), we computed the advection of the material lines up to the time \(t=2\). In the computation of advection, each material line is represented by 101 node points corresponding to \(s = i/100\times 0.1\,(i=0,\cdots, 100) \), and the advected node points are traced by using the Runge-Kutta method. Figure \ref{fig:torus_materialline} shows the position of the nodes of the three material lines at \(t=0, 1\) and \(2\). From this figure, we can see that each of the three material lines is elongated with deformation over time, but to evaluate the elongation in more detail, we also show the time evolution of the \textit{energy} \(\mathcal{E}\) of the material lines in figure \ref{fig:torus_materialline_energy}. Here, the energy \(\mathcal{E}\) of a material line \(\varphi(s;t)\,(s\in[0, s_0])\) is defined as: 
\begin{align}
    \mathcal{E}(t) := \frac{1}{2}\int_{0}^{s_0} \left|\frac{\partial \varphi}{\partial s}(\sigma; t)\right|^2 d\sigma.
    \label{k3}
\end{align}
Note that the word ``energy'' is not related to the kinematic energy of the flow, but it is a term of the differential geometry. Discretizing the right-hand side of \eqref{k3} by using \(I+1\) node points yields
\begin{align*}
    &\mathcal{E} (t) \approx \frac{1}{2} \sum_{i = 0}^{I-1} g(x(s_i), y(s_i))\left[\frac{(x(s_{i+1})-x(s_{i}))^2 + (y(s_{i+1})-y(s_{i}))^2}{(\Delta s)^2}\right] \Delta s, \\
    & s_{i} = i\Delta s\quad (i=0, \cdots, I), \qquad \Delta s = \frac{s_0}{I}.
\end{align*}
Therefore, for a fixed \(\Delta s\), the energy is proportional to the sum of the square lengths of the line elements into which the material line is discretized, and the square length of the line element is proportional to \(|\bm{\xi}|^2\) in \eqref{xinorm_evolution_M} and \eqref{xinorm_2nddiff_result}. For this reason, we use the energy of the material line for analysis rather than the length of the material line itself. We can see in figure \ref{fig:torus_materialline_energy} that the energy of each material line grows at an accelerated rate from \(t=0\) to \(t=1\). After that, the growth of the energy of the material lines slows down for \(x_0=0.30\times 2\pi\) and \(x_0 = 0.35\times 2\pi\), but for \(x_0 = 0.25\times 2\pi\), the growth of the energy remains accelerated. At \(t=2\), the material lines for \(x_0 = 0.30\times 2\pi\) and \(x_0 = 0.35\times 2\pi\) are located near the bend of the deformed vortex, indicating that the once-elongated vortex is reassembling. On the other hand, the material line for \(x_0 = 0.25\times 2\pi\) is located at the outer edge of the bent vortex and continues to be filamented. 

\section{Discussion and conclusions}\label{section_discussion}
In this paper, we have derived an expression for the strain acceleration tensor \(\mathsf{M}\), which describes the acceleration of the elongation of the tangent vector \(\bm{\xi}\) of the advected material line for a fluid governed by the Euler equation on a two-dimensional Riemannian manifold. To derive the expression, we have to compute the time derivative of \(\langle \mathsf{S}\bm{\xi}, \bm{\xi}\rangle\). This operation requires differential geometrical treatment because it involves spatial differentiation of the tangent vector \(\bm{\xi}\) of a material line or the velocity gradient tensor \(\nabla\bm{u}\) on a curved manifold. The resulting strain acceleration tensor \(\mathsf{M}\) contains a curvature term, which indicates that a material line is elongated or contracted due to the curvature of the manifold. The curvature effect works in the direction of promoting elongation of the material line if the curvature of the manifold is negative, and suppressing it if it is positive. Since the curvature tensor is intrinsic to Riemannian manifolds, such a term does not appear as long as we are dealing with flat planes. Note that we have only considered two-dimensional fluids, but we have not used any argument that depends on the dimension of the manifold in deriving \eqref{xinorm_2nddiff_result}. Therefore, the equation \eqref{xinorm_2nddiff_result} is valid for Euler flows on manifolds of arbitrary dimension. 

It should be noted that the appearance of the curvature term is due to the fact that we are considering a fluid subject to the Euler equation. The Euler equation on a Riemannian manifold is a geodesic equation for an infinite number of particles, corrected by adding a pressure term to satisfy incompressibility. We can therefore assume that the same curvature effects are at work when we consider the motion of a cluster of mass points. In fact, the form of the curvature term is consistent with the Jacobi equation. If we recall the derivation, this term arose from the covariant derivative of the advection term of the Euler equation \(\nabla_{\bm{u}}\bm{u}\) in the direction of \(\bm{\xi}\). Therefore, if the fluid follows an equation of motion with a similar advection term, a similar curvature term will appear in the expression of the strain acceleration. However, if the time evolution of the flow field is not governed by equations of motion with advection terms (such as the time-periodic flow field introduced in \citealp{pierrehumbert1991large}), there is no guarantee that the curvature effect will appear.

To determine the Lagrangian hyperbolic domain of Haller, it is necessary to check at each point on \(M\) whether \(\langle \mathsf{M}\bm{\xi},\bm{\xi}\rangle\) is positive definite on the set of vectors \(\bm{\xi}\) such that \(\langle \mathsf{S}\bm{\xi}, \bm{\xi}\rangle = 0\). Considering the result of the present study, the equation (5.1) of \cite{haller2005objective} cannot be used naively to determine Haller's hyperbolic domain on a non-flat manifold such as a sphere. There we have to consider the effect of curvature. Consider a point on a manifold where the sectional curvature is \(\kappa\), then the scale of the curvature radius of a curve passing the point is \(\kappa^{-1/2}\). Therefore, if the scale of the spatial change in the vorticity field is comparable to or larger than \(\kappa^{-1/2}\), the curvature effect can be significant in vortex dynamics. 

In section \ref{section_examples}, we analyzed three examples of vorticity fields using the equation \eqref{xinorm_2nddiff_result}. As the first example, we considered a vorticity field on a sphere consisting of two vorticity patches meeting on a latitude circle \(\mu = \mu_0\). This vorticity field represents a zonal jet whose axis is at \(\mu = \mu_0\). In this example, we can explicitly calculate the advection of the material line and its length. The time evolution of the length of the material line obtained explicitly in this way is consistent with the expression \(\frac{1}{2}d^2|\bm{\xi}|^2/dt^2\) obtained using \eqref{xinorm_2nddiff_result}. This is a very simple example that verifies the need to include a curvature term in the strain acceleration tensor. The vorticity field considered in this example is the same as that used by \cite{dritschel1988repeated}, where a weakly nonlinear theory was proposed to describe the initial stages of filamentation. The equations derived in the present study may help to quantify the subsequent filament elongation shown there. 

The second example, considered in subsection \ref{subsec_hypbsphere}, is a quadrupole state of stationary vortices on a sphere. Such a vorticity field, expressed as a superposition of spherical harmonics with total wavenumber 2, has the largest scale of all vorticity fields with zero angular momentum. This type of vorticity field is also known as a fixed energy solution of the minimum enstrophy theory \citep[e.g.,][]{herbert2013additional}. For this vorticity field, we computed Haller's hyperbolic domain using \eqref{xinorm_2nddiff_result} (top panel of figure \ref{fig:hypb_sphere}). Compared to Haller's hyperbolic domain obtained by neglecting the curvature term in \eqref{xinorm_2nddiff_result} (bottom panel of figure \ref{fig:hypb_sphere}), the true hyperbolic domain is slightly smaller. Since the sphere has a constant positive curvature everywhere, the curvature term in the strain acceleration tensor is negative semidefinite. Therefore, the hyperbolic domain with the curvature term properly considered is always smaller than the hyperbolic region with the curvature term ignored. To check whether the hyperbolic domain captures the Lagrangian coherent structure, we computed the finite-time Lyapunov exponent for the flow corresponding to the vorticity field \eqref{sphere_q_def} (figure \ref{fig:sphere_q}) and compared it with the time that each particle spent in the hyperbolic domain in the time interval (hyperbolicity time; figure \ref{fig:hypbtime}). The finite-time Lyapunov exponent is large in the two banded regions passing through \((\lambda,\mu)=(0,0)\) and \((\pi, 0)\), respectively. These regions are also indicated as the regions where the value of the hyperbolicity time is large. However, the hyperbolicity time also has annular maxima that are not seen in the distribution of the finite-time Lyapunov exponent. The hyperbolicity time for the hyperbolic domain calculated without the curvature term (bottom panel of figure \ref{fig:hypbtime}) shows that the annular maxima are even more pronounced. Small material lines containing particles starting from these annular maxima remain in the hyperbolic domain for a long time, but they continue to stretch and contract periodically. Therefore, the value of the Lyapunov exponent is not very large around the ring peaks. The periodic elongation and contracting of a material line that spends most of its time in the hyperbolic domain may seem to contradict the existence of a repelling/attracting material line, which \cite{haller2001lagrangian} and \cite{haller2005objective} demonstrated.  However, particles starting from these ring peaks pass through the cusp-like boundary of the hyperbolic region, and a small but finite length of material line partially extends outside the hyperbolic region. The existence of a repelling/attracting material line is not necessarily guaranteed for a time interval during which the cusp-like boundary is passed. In general, the hyperbolic domain with the curvature term ignored is slightly larger (smaller) for positive (negative) curvature than the hyperbolic domain with the curvature term properly included. Therefore, if there is a particle trajectory that stays at the edge of the correct hyperbolic domain for a long time, using a strain acceleration tensor that ignores the curvature term will lead to an overestimation (underestimation) of the value of the hyperbolicity time for positive (negative) curvature. 

The third example is the initial stage of time evolution from an initial vorticity field in a flow on a curved torus. We have used a curved torus for two reasons. The first is the ease of numerical computation. When the diagonal Riemannian metric is given, the equation relating the stream function and the vorticity is given by \eqref{torus_qlappsi}. Thus, the Laplace-Beltrami operator on the curved torus is the multiplication of \(g^{-1}\) and the Laplace-Beltrami operator on the flat torus. The same correspondence holds for the nonlinear term of the Euler equation. These correspondences make the numerical time integration easy, since the spectral method using trigonometric functions, which is used for the ordinary flat torus, is available. The other reason is the existence of points of negative curvature. At points of negative curvature, the curvature term in \eqref{xinorm_2nddiff_result} is positive, and thus the negative curvature may contribute to the elongation of the vortices. In the third example, a positive vortex at the initial time was deformed in the negative curvature region, and part of it became filamentary. By comparing the values of the terms in \eqref{xinorm_2nddiff_result} it was shown that the filamentation was induced by the negative curvature of the region. We also calculated the advection of the material line and the time evolution of their energy. Note that in the right half of the positive vortex (\((x,y)=(0.35\times 2\pi, 0.70\times 2\pi)\)), the pressure term of \eqref{xinorm_2nddiff_result} is positive. The elongation of material lines is due not only to the curvature effect but also to the spatial variation of the pressure gradient forces. The material lines starting from the right half of the vortex and its center (near \((x, y) = (0.30 \times 2\pi, 0.70 \times 2\pi)\)) then enter the bent region of the deformed vortex, where the curvature is close to zero and the bent vortex begins to coalesce, and the elongation of the material lines slows down (figure \ref{fig:torus_materialline_energy}). On the other hand, in the left part of the vortex (around \((x, y) = (0.25 \times 2\pi, 0.70 \times 2\pi)\)), the strain acceleration is positive only if the curvature term is taken into account. The elongation of the material line from this part of the vortex is accelerated at least up to time \(t = 2\). This again shows the importance of surface curvature. 

\begin{figure}
    \centering
    \includegraphics[width=0.5\linewidth]{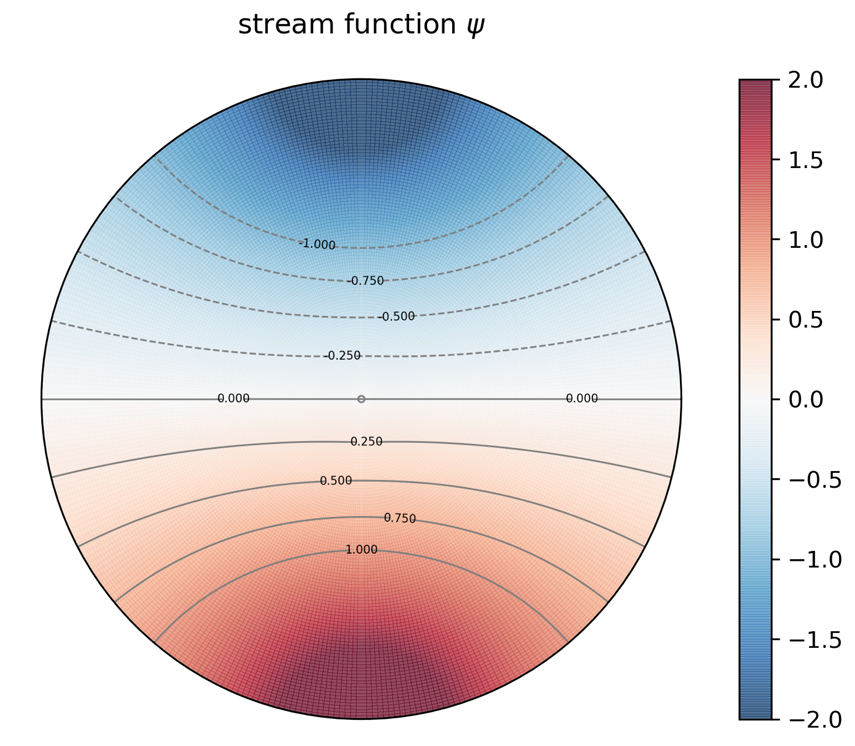}
    \caption{The stream function \(\psi\) of the flow defined by \eqref{flow_on_pd}. The streamlines corresponding to \(\psi = -1.0, -0.75, -0.5, -0.25, 0, 0.25, 0.50, 0.75\) and \(1.0\) are also shown.}
    \label{fig:pdisk_stf}
\end{figure}

In the third example, the deformation of the vortex was caused not only by the flow induced by the vortex itself but also by the flows induced by other nearby vortices. At \(t=0\), the positive vortex to be deformed is located between two vortices of opposite sign, the positive vortex located near \((x,y) = (0.95\times 2\pi, 0.55\times 2\pi)\) and the negative vortex located near \((x,y) = (0.60\times 2\pi, 0.60\times 2\pi)\) (see figures \ref{fig:torus_initvol} and \ref{fig:torus_init_press_okuboweiss}). 
These two vortices induce a jet in the positive direction of the \(y\) axis. The filamentation of the vortex occurs during advection by this jet. Of course, even on a flat plane, the streamlines of a jet formed by two vortices of different signs will spread out toward the outlet, and the reader might think that the vortex deformation is not due to the negative curvature. However, due to the negative curvature, the jet streamlines have a wider spread. To understand this, consider a steady flow on the Poincar\'{e} disk \(D\) as a toy model. Recall that the Poincar\'{e} disk is a unit open disk in the complex plane, \(D = \{x + \sqrt{-1}y\in \mathbf{C}\mid x,y \in \mathbf{R},\, x^2+y^2<1\}\), endowed with the Riemannian metric \(ds^2 = 4(1-|z|^2)^{-2}|dz|^2\). The Poincar\'{e} disk is a Riemannian manifold with a constant negative curvature \(R_{1221} = -1\), and the geodesics are arcs orthogonal to the boundary circle \(\partial D\) at two points. The Euler equation on \(D\) has the same form as \eqref{torus_eulereq_u}--\eqref{torus_non_divergent} if we set \(g = 4/(1-|z|^2)^2\). The flow field defined below is an incompressible and irrotational flow on \(D\):
\begin{align}
    u(x,y) - \sqrt{-1}v(x,y) = \frac{1-|z|^2}{2}f(z);\qquad f(z)=\frac{2}{z^2+1}.\label{flow_on_pd}
\end{align}
Note that the flow field becomes incompressible and irroational if and only if the function \(\sqrt{g}(u-\sqrt{-1}v)\) is a holomorphic function of \(z\). This flow field represents a jet in the positive direction of the \(x\) axis and the streamlines match the geodesics on \(D\) (figure \ref{fig:pdisk_stf}). Given this flow as a background flow, a weak vortex blob placed at \(z=0\) will be advected in the positive direction of the \(x\) axis with rapid elongation in the \(y\) axis direction. Note that \(\mathsf{H}(p)\) is positive semidefinite and \(\nabla\bm{u} = 0\) at \(z=0\). The function \(f(z)\) has two poles at \(z = \pm \sqrt{-1}\), but they are at ``points at infinity'' and the flow domain \(D\) has no singular points. This flow field gives an example of flows such that the fluid particles move along geodesics, but their trajectories spread out. The existence of such flows is a characteristic of negative curvature manifolds; manifolds with positive curvature or flat planes do not allow such flows without singular points. In other words, a unidirectional jet on a negative curvature manifold spreads out due to the curvature of the flow domain. Of course, the flow considered in section \ref{subsec_torus} is not as simple as this toy model flow, but we can confirm the contribution of the negative curvature to the acceleration of the vortex elongation from the analysis of figure \ref{fig:torus_hypb_term}.

Before concluding this paper, let us discuss the significance of this entire study and its prospects in the following paragraphs. Equation \eqref{xinorm_2nddiff_result} extends the Lagrangian coherent structure criterion proposed by Haller to fluids on curved surfaces. The analysis of the Lagrangian coherent structure is not only a classical problem but also an important problem in the real world. For example, tracer transport by the polar vortex in the Earth's stratosphere has been studied \citep[e.g.,][]{mizuta2001chaotic}. In such studies, it is important to distinguish where mixing is likely to occur. For vortices of planetary scale such as the polar vortex, the curvature of the planet may have some effects to the transport of the vorticity or materials. The results of this study may be applicable to such studies. 

In addition, as an impact on more fundamental research, Equation \eqref{xinorm_2nddiff_result} can motivate the study of fluid dynamics on manifolds of negative curvature, beyond the usual framework of fluids on flat or spherical surfaces. As shown in the third example of section \ref{section_examples}, at points of negative curvature, the elongation of the material line is accelerated by the geometric effect. In two-dimensional fluids, the process of filamentation of the vorticity field is an important process of energy transfer to larger scales, and it may be interesting to investigate whether this energy transfer process is accelerated by the curvature of the domain. 

The hyperbolic Lagrangian coherent structure can also lead to chaotic advection. As a geometric theory of chaos in flow fields, Arnold's theory of sectional curvature of diffeomorphism groups \citep{arnold2009topological} is famous. In this context, ``sectional curvature'' refers to the curvature of the infinite-dimensional Lie group of a volume-preserving diffeomorphism on a manifold, not to the curvature of the domain itself. In the case of the flat torus, the relationship between the curvature of the infinite dimensional Lie group and the elongation of the material lines has been studied \citep{hattori1994motion}. The relationship between the curvature of the domain itself and the sectional curvature of the diffeomorphism group may be a subject for future research. 

The motion of a fluid on a manifold with negative curvature may also be an interesting future research topic in the context of the pattern formation of a two-dimensional fluid. 
It has been pointed out that the motion of a two-dimensional fluid is not ergodic, and even turbulent flows do not necessarily reach statistical equilibrium, defined as the maximizer of the mixing entropy, on the flat torus and spherical surfaces \citep{segre1998late,morita2011collective,ryono2024statistical}. However, it is not clear whether a statistical equilibrium solution is actually realized when mixing is sufficiently promoted. To address this issue, the following strategy can be proposed, taking into account the results of this study. If we allow the Riemannian metric to vary with time, it can create a region of negative curvature, which can promote mixing. In particular, if the Riemannian metric is varied so that the volume element is preserved, then the energy and all Casimir invariants of the flow are conserved. For example, we can introduce a time-periodic Riemannian metric such as:
\begin{align*}
     g = \left(
    \begin{array}{cc}
        \gamma(x,y;t) & 0 \\[1.5em]
        0 & \gamma(x,y;t)^{-1}
    \end{array}
    \right).
\end{align*}
In this way, we can give some change to the motion of the vortices without violating the symplectic structure of the Euler equation, and this change could promote the mixing process. Note, however, that the formulas \eqref{xinorm_evolution_M} and \eqref{xinorm_2nddiff_result} should be modified in this case to include the effects of the time change in the metric.

In conclusion, the results of this paper are not just an extension of fluid dynamics on a plane to Riemannian manifolds, but emphasize that the Euler equations have their origin in the geodesic equations, which gives a new perspective on the relationship between the geometry of the flow domain and the motion of the flow field in terms of the curvature of the domain and the elongation of the vortices.

\appendix
\section{The derivation of \eqref{xi_evolution_M}}
By the flow field \(\bm{u}\), a small neighborhood \(U\subset M\) of the point \(x_0 = \varphi(0; 0)\) is advected homeomorphically. Thus, if a particle starting from \(x\in U\) at initial time reaches \(\Phi_t (x)\) at time \(t\), the map \(\Phi_t\) is a diffeomorphism. Considering the invariance of the Euler equation under the time reversal transformation, the map \(\Phi_t\) can also be defined for \(t<0\).

We can take a local coordinate \((s_1,\cdots, s_n)\) of \(U\) such that the parameter \(s\) of the material line corresponds to \(s_1\), and the point \(x_0\) is expressed as \((s_1,\cdots,s_n)=(0,\cdots,0)\). Furthermore, the same coordinate can be used for \(\Phi_t(U)\) at any \(t\) since \(\Phi_t\) is a diffeomorphism. The vector field \(\bm{\xi}(t)\) defined on the trajectory of the particle \(t\mapsto \varphi(0;t)\) can be expressed as:
\begin{align*}
    \bm{\xi}(t) = \left(\frac{\partial}{\partial s_1}\right)_{\varphi(0;t)}=(D\Phi_t)\left(\frac{\partial}{\partial s_1}\right)_{\varphi(0;0)}
\end{align*}
by definition. Here, \(D\Phi_t\) is the differential map of \(\Phi_t\). As we did in section \ref{subsubsec_stretch}, we consider the extended space \(E=M\times (-\tau,\tau)\). The set \(U_E =  \{(\Phi_t(x),t)\in E\mid x\in U, -\tau<t<\tau\}\) is an open submanifold of \(E\) and has \((s_1,\cdots,s_n,t)\) as its coordinate. Consider the vector field on \(U_E\) defined as:
\begin{align*}
    &\tilde{\bm{u}}(\Phi_t(x), t) = \bm{u}(\Phi_t(x),t) + \left(\frac{\partial}{\partial t}\right)_{(\Phi_t(x),t)}\\
    &\tilde{\bm{\xi}}(\Phi_t(x),t) = \left(\frac{\partial}{\partial s_1}\right)_{\Phi_t(x)}.
\end{align*}
Since the coordinate \(s_1\) is frozen to particles, the two diffeomorphisms generated by \(\tilde{\bm{u}}\) and \(\tilde{\bm{\xi}}\) commute. Therefore the vector fields \(\tilde{\bm{u}}\) and \(\tilde{\bm{\xi}}\) also commute as follows:
\begin{align*}
   [\tilde{\bm{u}},\tilde{\bm{\xi}}] = 0.    
\end{align*}
By the symmetry of the Levi-Civita connection, \([\tilde{\bm{u}},\tilde{\bm{\xi}}] = \nabla^E_{\tilde{\bm{u}}}\tilde{\bm{\xi}} - \nabla^E_{\tilde{\bm{\xi}}}\tilde{\bm{u}}\). By setting \(s_1 = \cdots = s_n = 0\), we obtain
\begin{align*}
    \frac{D\bm{\xi}}{dt} &= \nabla_{\tilde{\bm{u}}}\tilde{\bm{\xi}}|_{(s_1,\cdots,s_n,t)=(0,\cdots,0,t)} \\
    &=\nabla_{\tilde{\bm{\xi}}}\tilde{\bm{u}}|_{(s_1,\cdots,s_n,t)=(0,\cdots,0,t)} =\nabla_{\bm{\xi}}\bm{u}.
\end{align*}

\section{Evaluation of the spectral coefficients of the Riemannian metric}
Here we will show that the evaluation \eqref{metric_g_evaluation} holds for the following metric and its expansion coefficients:
\begin{align*}
    g(x,y) = e^{\alpha\sin x\sin y};\qquad a_{k,l}=\frac{1}{(2\pi)^2}\int_{0}^{2\pi}\int_{0}^{2\pi}g(x,y)e^{-\sqrt{-1}(kx+ly)}dxdy.
\end{align*}
First, we note that
\begin{align*}
    a_{k,l} &= \frac{1}{(2\pi)^2}\int_0^{2\pi}\int_0^{2\pi} \sum_{n=0}^\infty \frac{(\alpha \sin x\sin y)^n}{n!}e^{-\sqrt{-1}(kx+ly)}dxdy \\ 
    &= \frac{1}{(2\pi)^2}\sum_{n=0}^\infty\frac{\alpha^n}{n!}\int_0^{2\pi}\int_0^{2\pi}\sin^n x\sin^n y\,e^{-\sqrt{-1}(kx+ly)}dxdy \\ 
    &= \sum_{n=0}^\infty \frac{\alpha^n}{n!} I_{n,k} I_{n,l},
\end{align*}
where 
\begin{align*}
    I_{n,k} = \frac{1}{2\pi}\int_{0}^{2\pi} \sin^n x\,e^{-\sqrt{-1}kx} dx.
\end{align*}
Using Euler's formula, we obtain
\begin{align*}
    I_{n,k} &= \frac{1}{2\pi} \int_{0}^{2\pi} \left(\frac{e^{\sqrt{-1}x}-e^{-\sqrt{-1}x}}{2\sqrt{-1}}\right)^n e^{-\sqrt{-1}kx} dx \\
    &=\frac{1}{(2\sqrt{-1})^n}\sum_{m=0}^n\binom{n}{m} \frac{1}{2\pi}\int_{0}^{2\pi}e^{\sqrt{-1}mx}e^{-\sqrt{-1}(n-m)x}e^{-\sqrt{-1}kx}dx\\
    &=\frac{1}{(2\sqrt{-1})^n}\sum_{m=0}^n\binom{n}{m} \frac{1}{2\pi}\int_{0}^{2\pi} e^{\sqrt{-1}(2m-n-k)x}dx.
\end{align*}
The integration on the rightmost-hand side has a non-zero value only for \(m=0,\cdots, n\) such that \(2m-n-k=0\). Therefore, if \(n+k\) is odd or \(n+k<0\), then \(I_{n,k}=0\). We thus consider the case that \(n+k\) is even and non-negative. Note that \(2m-n-k=0\) if and only if \(m=(n+k)/2\) and \((n+k)/2\leq n\) only if \(k\leq n\). The value \(I_{n,k}\) is nonzero only if \(|k|\leq n\) and \(k+n\) is even, and then
\begin{align*}
     I_{n,k} = \frac{1}{(2\sqrt{-1})^n}\binom{n}{\frac{n+k}{2}}.
\end{align*}
Therefore, 
\begin{align*}
    a_{k,l} = \sum_{n}{}' \frac{\alpha^n}{n!} \frac{1}{(2\sqrt{-1})^{2n}}\binom{n}{\frac{n+k}{2}}\binom{n}{\frac{n+l}{2}},
\end{align*}
where \(\sum'\) means the summing over the integers \(n\geq \max \{|k|, |l|\}\) such that both \(k+n\) and \(l+n\) are even. In particular, if \(k\) and \(l\) have different parities, \(a_{k,l}\) is zero. When \(k\) and \(l\) have the same parity, 
\begin{align*}
    a_{k,l} = \sum_{j=0}^\infty b_j,\qquad b_j = \frac{\alpha^{L+2j}}{(L+2j)!}\frac{1}{(2\sqrt{-1})^{2(L+2j)}} \binom{L+2j}{\frac{L+2j+k}{2}}\binom{L+2j}{\frac{L+2j+l}{2}}
\end{align*}
for \(L=\max\{|k|,|l|\}\). Since for \(j=0,1,2,\cdots\) the inequality
\begin{align*}
    \frac{|b_{j+1}|}{|b_{j}|} &= \frac{\alpha^{2}}{4^2(L+2j+2)}\frac{L+2j+2}{\frac{L+2j+1}{2}+1} \\
    &=\frac{\alpha^2}{8(L+2j+3)}\leq \frac{\alpha^2}{8(L+3)}
\end{align*}
holds, for sufficiently large \(k\) and \(l\), we have
\begin{align*}
    |a_{k,l}|&\leq \sum_{j=0}^\infty |b_j| \leq |b_0| \frac{1}{1-\frac{\alpha^2}{8(L+3)}}\\
    &\leq 2|b_0| = \frac{2}{L!}\left(\frac{\alpha}{4}\right)^{L}\binom{L}{\frac{L+k}{2}}\binom{L}{\frac{L+l}{2}}.
\end{align*}
Here, \(L\) is equal to \(|k|\) or \(|l|\). If \(L=|k|\), then we obtain
\begin{align*}
    |a_{k,l}|\leq \frac{2}{L!}\left(\frac{\alpha}{4}\right)^{L} \binom{L}{\frac{L+l}{2}} \leq \frac{2}{L!}\left(\frac{\alpha}{4}\right)^L \binom{L}{\frac{L}{2}} \leq C\frac{1}{L!}\left(\frac{\alpha}{2}\right)^L,
\end{align*}
where \(C>0\) is a constant that does not depend on \(k\) and \(l\). The same evaluation holds for the case of \(L=|l|\). Note that the last inequality follows from Wallis's formula:
\begin{align*}
     \binom{2n}{n} \approx \frac{4^n}{\sqrt{n}}\leq 4^n \qquad (n\to\infty).
\end{align*}
We can derive the same form of the inequality for the Fourier coefficients of \(g^{-1}\) in a similar way.

\section{Explicit specification of the initial vorticity field on the curved torus}
In section \ref{subsec_torus}, we have considered an example of the flow on a curved torus, which represents the onset of turbulence. Although the initial vorticity field considered was originally generated by random numbers, we explicitly write down its spectral coefficients here for the sake of reproducibility.

The initial vorticity field is given by the following equation:
\begin{align*}
    q(x,y) = \sum_{k=0}^2\sum_{l=0}^2 &(a_{k,l}\cos kx \cos ly + b_{k,l} \cos kx\sin ly \\
    &\qquad+ c_{k,l} \sin kx \cos ly + d_{k,l}\sin kx \sin ly),
\end{align*}
where 
\begin{align*}
\begin{array}{ll}
a_{0,0} = 0.05983618385516437 & a_{0,1} = -0.07844948104163184 \\
a_{0,2} = 0.12523658449465810 & a_{1,0} = 0.16532280170939714 \\
a_{1,1} = -0.11693416269523071 & a_{1,2} = 0.39389003962238250 \\
a_{2,0} = -0.08428800796409970 & a_{2,1} = 0.04229008794099242 \\
a_{2,2} = 0.30077848482043962 & \qquad \\ 
\qquad & \qquad \\
b_{0,1} = -0.19049291872705523 & b_{0,2} = 0.05827967227392698 \\
b_{1,1} = -0.10944106134331011 & b_{1,2} = 0.32536644789854924 \\
b_{2,1} = 0.26029577962448952 & b_{2,2} = 0.00188425339698650 \\
\qquad & \qquad \\ 
c_{1,0} = 0.19681063366118581 & c_{1,1} = -0.31923881592497894 \\
c_{1,2} = -0.16743435702548920 & c_{2,0} = -0.10316648330130607 \\
c_{2,1} = 0.02821210104193395 & c_{2,2} = -0.22662753731338633 \\
\qquad & \qquad \\ 
d_{1,1} = -0.19180157757767999 & d_{1,2} = 0.23382559599065009 \\
d_{2,1} = 0.39408167049915921 & d_{2,2} = -0.37903947388947745. \\
\end{array}    
\end{align*}

\section*{Acknowledgment}
This work was supported by JSPS KAKENHI, Grant Number JP24KJ1340.

\bibliographystyle{jphysicsB}
\bibliography{ref}

\end{document}